%% file: 102153.tex
\magnification=\magstephalf

\input 102153macro

\overfullrule=0pt

\input 102153ref


\afont
\noindent
Convergent Asymptotic Expansions of Charlier,
\noindent
Laguerre and Jacobi Polynomials

\bigskip

{\leftskip=25pt
\afont
\noindent
Jos\'e L. L\'opez$^1$  and Nico M. Temme$^2$

\medskip\noindent
\rm
$^1$
Departamento de Mat\'ematica e Inform\'atica,
Universidad P\'ublica de Navarra,
31006-Pamplona,
Spain ({\tt jl.lopez@unavarra.es}),

\noindent
$^2$
CWI,
P.O. Box 94079,
1090 GB Amsterdam,
The Netherlands ({ \tt nicot@cwi.nl}).
\medskip

\bfont
\bigskip
\noindent
Convergent expansions are derived for three types of orthogonal
polynomials: Charlier, Laguerre and Jacobi. The expansions have asymptotic
properties for large values of the degree. The expansions are given in
terms of functions that are special cases of the given polynomials.
The method is based on expanding integrals in one or two points of the
complex plane, these points being saddle points of the phase functions
of the integrands.

\vskip 0.6cm \noindent
\cfont
2000 Mathematics Subject Classification:
\bfont
33C45, 41A60, 30E20.
\par\noindent
\cfont
Keywords \& Phrases:
\bfont
Charlier polynomials,
Laguerre polynomials,
Jacobi polynomials,
asymptotic expansions,
saddle point methods,
two-points Taylor expansions.
\par\noindent
\cfont
}
\rm
\parindent=15pt
%
%
\newcap{1}{Introduction}

\noindent
In a previous paper \temmelopez\ we have studied the expansion of an
analytic function at two finite points in the complex plane.  The domain of
convergence is a Cassini oval around the two points.  The main motivation
for that paper was to obtain the coefficients of asymptotic expansions of
certain integrals.  In the present paper we give a few examples in which
the expansion of an integral at two saddle points yields a convergent
expansion that has an asymptotic property for large values of a parameter.

In the well-known methods for deriving asymptotic expansions of integrals a
basic step is transforming the integral into a standard form, and the
transformation usually gives a new integral in which the integrand contains
implicitly defined functions that are difficult to handle.  In the method
of this paper we avoid a transformation and, in addition, we derive
convergent expansions.

We start with a simple example in which only one saddle point occurs, and
in which a function is expanded at that saddle point.  This gives an
expansion for the Charlier polynomials.

In two other examples (Laguerre and Jacobi polynomials) we take into
account two saddle points, and again two convergent expansions can be
constructed with the desired property.  The approximants belong to the same
class of polynomials as the original ones, but they are of a simpler type
(Hermite and Chebyshev, respectively).  The asymptotic property follows
from recursion relations for functions appearing in the expansions.  The
convergence follows from the fact that an integral along a finite contour
is expanded inside a domain of uniform convergence.

In the examples given in this paper the contour integrals are based on
Cauchy-type integrals obtained from generating functions. When the
contour is finite a proof of the convergence is usually rather easy.
For more general finite contours and more general integrals we expect
that the method can be applied as well.  For example, we can apply
the method to the Gauss hypergeometric function and the incomplete
gamma function, with different integral representations.

Also, the methods of this paper can be generalized by considering
Taylor expansions at more than two points. In  \multayl\ we give
details on the theory of multi-point Taylor expansions, and in a
future paper we will give details on applications to integrals with,
for example, three saddle points.

We show a few graphs that indicate the nature of the approximations,
and in a final section we mention a few examples in which other functions
are considered.

\newcap{2}{A simplified version of the saddle point method}

\noindent
Throughout this paper we are concerned with finding asymptotic expansions of
integrals of the form
\eqn\inte{
F(n)\equiv\int_\Gamma f(w)e^{ng(w)}{dw\over w^{n+1}},
}
where $f(w)$ and $g(w)$ are analytic in a domain $\Omega$ of the complex plane
that contains the origin; $\Gamma$ is a circle with
center at the
origin and contained in $\Omega$; $n$ is a large positive integer.
We assume, as it usually happens to be the case, that the asymptotic
behavior of the integral $F(n)$ for large $n$ is determined by
contributions from the saddle points of $\varphi(w)=g(w)-\ln w$
$[$\wong, Chap. 2, \S 4$]$. 

The standard saddle point method consists of
{\parindent=25pt
\item{(i)\sq}
      deforming the contour of integration
      $\Gamma$ into a new path that crosses one or some of the saddle points
      of $\varphi(w)$;
\item{(ii)\quad}
       a suitable change of the variable of integration;
\item{(iii)\quad}
       application of Watson's lemma or Laplace's method.
\par
}

Instead of applying the
standard saddle point method,
we will proceed in a simpler way: just substitute a power series expansion
at one or more saddle points of the function $f(w)$
in \inte. If there is just one saddle point $w_0$, then that
power series is its Taylor
expansion at $w_0$:
\eqn\simple{
f(w)=\sum_{k=0}^\infty{f^{(k)}(w_0)\over k!}(w-w_0)^k,
}
which is uniformly convergent for $w$ in a disk
$D_r(w_0)\equiv\lbrace w\in\Omega$, $\vert w-w_0\vert<r\rbrace$
with center at $w_0$ and radius
$r=$ Inf$_{w\in\Cs\setminus\Omega}\vert w-w_0\vert$.
If there are two saddle points $w_1$ and $w_2$, then that
power series is its two-point Taylor series at $w_1$ and $w_2$ \temmelopez:
\eqn\double{
f(w)=\sum_{n=0}^\infty\left\lbrack a_n(w-w_1)+a_n'(w-w_2)\right\rbrack
(w-w_1)^n(w-w_2)^n,
}
where
\eqn\ao{
a_0\equiv{f(w_2)\over w_2-w_1}, \hskip 2cm a_0'\equiv{f(w_1)\over w_1-w_2}
}
and,  for $n=1,2,3,...$,
\eqn\coefibis{
a_n\equiv{1\over n!}\sum_{k=0}^n{(n+k-1)!\over
k!(n-k)!}{(-1)^{n+1}nf^{(n-k)}(w_2)+(-1)^kkf^{(n-k)}(w_1)\over
(w_1-w_2)^{n+k+1}},
}
\eqn\aprime{
a_n'\equiv{1\over n!}\sum_{k=0}^n{(n+k-1)!\over
k!(n-k)!}{(-1)^{n+1}nf^{(n-k)}(w_1)+(-1)^kkf^{(n-k)}(w_2)\over
(w_2-w_1)^{n+k+1}}.
}
The expansion \double\ is uniformly convergent for $w$ in a Cassini oval
$$
O_r(w_1,w_2)\equiv\lbrace w\in\Omega,\ \vert w-w_1\vert\vert 
w-w_2\vert<r\rbrace
$$
with foci at $w_1$ and $w_2$ and "radius"
$r={\rm Inf}_{w\in\Cs\setminus\Omega}\lbrace\vert w-w_1\vert\vert
w-w_2\vert\rbrace$; see  \temmelopez.

If we substitute now \simple\ or \double\ in \inte\ and interchange
summation and integration we obtain
an expansion of $F(n)$. This is proved in the following two propositions.

\noindent
{\hfont Proposition 2.1.} {\it Let  the right-hand side of \simple\ converge
uniformly to $f(w)$ for $w\in
D_r(w_0)$ with
$\vert w_0\vert<r$, then
\eqn\uno{
F(n)=\sum_{k=0}^\infty{f^{(k)}(w_0)\over
k!}\int_\Gamma(w-w_0)^ke^{n\varphi(w)}{dw\over w}.
}
}

\noindent
{\it Proof.} If $\vert w_0\vert<r$, then $0\in D_r(w_0)$. Then we can
choose a small
enough circle $\Gamma$ in \inte\ such that $\Gamma\in D_r(w_0)$. Therefore,
expansion \simple\
is uniformly convergent for $w\in\Gamma$. Introducing \simple\ in \inte\
and interchanging summation and integration we obtain \uno.
\hfill$\boxe$

\noindent
{\hfont Proposition 2.2.} {\it Let the right-hand side of \double\ converge
uniformly to $f(w)$ for $w\in
O_r(w_1,w_2)$ with $\vert w_1w_2\vert<r$, then
\eqn\dos{\eqalign{
F(n)=&\sum_{k=0}^\infty
a_k\int_\Gamma(w-w_1)^{k+1}(w-w_2)^ke^{n\varphi(w)}{dw\over
w}+ \cr &
\sum_{k=0}^\infty 
a_k'\int_\Gamma(w-w_1)^{k}(w-w_2)^{k+1}e^{n\varphi(w)}{dw\over
w}. \cr} }
}

\noindent
{\it Proof.} The proof is similar to that of Proposition 2.1.
\hfill$\boxe$

In the remaining part of the paper we apply Proposition 2.1 or 
Proposition 2.2
to three specific examples of integrals $F(n)$ representing
Charlier polynomials $C_n^a(nx)$,
Laguerre polynomials $L_n^\alpha(nx)$ and Jacobi polynomials
$P_n^{(\alpha,\beta)}(x)$. In this way,
we obtain expansions of these polynomials for large values of $n$.
In each example separately we prove that the corresponding expansions
\uno\ or \dos\
are convergent in a certain region of the variable $x$ and that
in fact, they have an
asymptotic nature for large $n$, uniformly with respect to $x$ in certain
domains of the
region of convergence.

\newcap{3}{Asymptotic expansions of Charlier polynomials
in terms of Gamma functions}

\noindent
The Charlier polynomials are defined by the generating function
\eqn\gfchar{
e^{-aw}(1+w)^{x}=\sum_{n=0}^\infty C_n^a(x){{w^n}\over{n!}},
}
and have the explicit expression
\eqn\explchar{
C_n^a(x)=\sum_{k=0}^n{n\choose k}{x\choose k} k! (-a)^{n-k}.
}
The Charlier polynomials are orthogonal with respect to a
discrete distribution on the positive real line.
For an overview of properties, see \koek. Recent papers on asymptotics
are \bowong, \dunster, \goh, and \hsu.
In this section we give a simple convergent expansion of $C_n^a(xn)$,
that has an asymptotic property for large $n$, uniformly for complex
$a$ in compact sets and for complex $x$  bounded away from 1.
A uniform expansion that holds for $x$ in a compact
neighbourhood of $x=1$ is given in \bowong, where
$J-$Bessel functions are used in the approximations.
Also, uniform expansions that hold for $-\infty<x<\infty$ are
given in \dunster.

\noindent
{\hfont Theorem 3.1.} {\it For $x\ne 1$, $a\in\CC$ and $n\in\Ns$, the Charlier
polynomials have the
expansion
\eqn\expanchar{
C_n^a(xn)=e^{a/(1-x)}\sum_{k=0}^\infty{(-a)^k\over k!}\Phi_k(x,n),
}
where
\eqn\fio{
\Phi_0(x,n)\equiv{\Gamma(nx+1)\over\Gamma(nx+1-n)}, \hskip 1cm
\Phi_1(x,n)\equiv {\Phi_0(x,n)\over(1-x)(n(x-1)+1)}
}
and, for $k=0,1,2,...$,
\eqn\fichar{
\Phi_k(x,n)\equiv{\Gamma(nx+1)\over\Gamma(nx-n+1)}{1\over(1-x)^k}
{}_2F_1(-k,-n,nx-n+1;1-x).
}
The sequence $\lbrace\Phi_k(x,n)$, $k=0,1,2,...\rbrace$ satisfies the 
recurrence
\eqn\recuchar{
\Phi_k(x,n)={1\over n(x-1)+k}\left\lbrack{x(1-k)-k\over x-1}\Phi_{k-1}(x,n)+
{x(1-k)\over(x-1)^2}\Phi_{k-2}(x,n)\right\rbrack,
}
where $k=2,3,\ldots \,.$ and is an asymptotic sequence for large $n$. For
fixed $k$ we have
\eqn\fiasy{
\Phi_k(x,n)={\cal O}(n^{-\lfloor (k+1)/2\rfloor}(n\vert x\vert)_n)
}
when $n\to\infty$,
where $\lfloor \alpha \rfloor$ is the integer part of the real number
$\alpha$. The asymptotic property holds uniformly with respect to complex
$x$, $|x-1|\ge\delta>0$.
}

Observe that in the notation of Pochhammer's symbol $(a)_n$, defined by
\eqn\poch{
(a)_0=1, \quad (a)_n=\frac{\Gamma(a+n)}{\Gamma(a)}= a(a+1)\cdots (a+n-1), 
\quad n=0,1,2,\ldots,
}
we have
\eqn\fifi{
\Phi_0(x)=(-1)^n{\Gamma(n-nx)\over \Gamma(-nx)}=
(-1)^n (-nx)_n=(nx)(nx-1)\cdots(nx-(n-1)).
}

\noindent
{\it Proof.} From \gfchar\ we derive the integral
representation
$$
C_n^a(x)={n!\over 2\pi i}\int_\Gamma e^{-aw}(1+w)^x{dw\over w^{n+1}},
$$
where $\Gamma$ is a circle with center at the origin and radius $<1$.
We write this in the form
\eqn\charly{
C_n^a(nx)={n!\over 2\pi i}\int_\Gamma e^{-aw}e^{n\varphi(x,w)}{dw\over w},
}
where $\varphi(x,w)\equiv x\log(1+w)-\log w$. The only saddle point of
$\varphi(x,w)$ is
$w_0=(x-1)^{-1}$. The function $e^{-aw}$ is an entire
function
of $w$. Hence, the expansion
\eqn\expo{
e^{-aw}=\sum_{k=0}^\infty{(-a)^ke^{-aw_0}\over k!}(w-w_0)^k
}
is locally uniformly convergent for $w\in \Cs$ and $x\ne 1$.
Therefore, after substituting this expansion in
\charly\ and using Proposition 2.1 we obtain \expanchar\ with
\eqn\ficharly{
\Phi_k(x,n)={n!\over 2\pi i}\int_\Gamma(w-w_0)^k(1+w)^{xn}{dw\over w^{n+1}}.
}
To obtain the recurrence \recuchar\ we write
$$
\Phi_k(x,n)={(n-1)!\over 2\pi i}{1\over x-1}\int_\Gamma(w-w_0)^{k-1}
(1+w){\partial e^{n\varphi(x,w)}\over\partial w}\,dw.
$$
Integrating by parts and performing a few straightforward manipulations we
obtain \recuchar. Equalities \fio-\fichar\
follow after simple calculations. The asymptotic behavior in \fiasy\
for large $n$ follows from \fio\ and \recuchar.
 From \fio\ we see that $\Phi_0={\cal O}((n\vert x\vert)_n)$ (see formula
(19) below) and that
  $\Phi_1={\cal O}((n\vert x\vert)_n/n)$. Therefore,
  \fiasy\ is true for $k=0,1$. From here, the proof
follows by induction over $k$. If \fiasy\ holds up to $k$, then
$\Phi_{k-1}={\cal O}(n^{-[k/2]}(n\vert x\vert)_n)$ and
$\Phi_k={\cal O}(n^{-[(k+1)/2]}(n\vert x\vert)_n)$.
Using this in \recuchar\ with $k$
replaced by $k+1$ we have that
$\Phi_{k+1}={\cal O}(n^{-[(k+2)/2]}(n\vert x\vert)_n)$.
\hfill$\boxe$

Property \fiasy\ holds
uniformly for $|x-1|\ge\delta>0$, also for complex $x$. Also, detailed
information on the
asymptotic behavior can easily be obtained from \ficharly.
\hfill $\boxe$

From \fifi\ it follows that $\Phi_0(x)$ has $n$ zeros at
$x=m/n$, $m=0,1, \ldots n-1$ ($\Phi_1(x)$ has the same zeros,
except for $x=(n-1)/n$).
See also Table 3.2, where we give
numerical values of the zeros of $C_n^a(nx)$ for $n=10$ and $a=1$.
 From the graphs in Figure 3.3 we also see that
the early zeros are approximated quite well.
\topinsert
$$\vbox{\hsize=10truecm\offinterlineskip
   \vskip3pt
$$\vbox{
\halign{\strut
\quad
  $#$\quad &  \qquad $#$ \qquad\hfil &   \qquad $#$ \qquad \hfil  \cr
     &\  \, 0.000000090   & \ \,  0.534449998 \cr
     & \ \, 0.100006223   & \ \,  0.680932968 \cr
     & \ \, 0.200157621   & \ \,  0.855641877 \cr
     & \ \, 0.301812498   & \ \,  1.068772397 \cr
     & \ \,  0.410358953  & \ \,  1.347867376 \cr
}}$$
\centerline{{\hfont Table 3.2.}\sq The zeros of $C_n^a(nx)$ for $n=10$ and
$a=1$.}
}$$

\endinsert

\bigskip
\centerline{{\epsfxsize=6cm \epsfbox{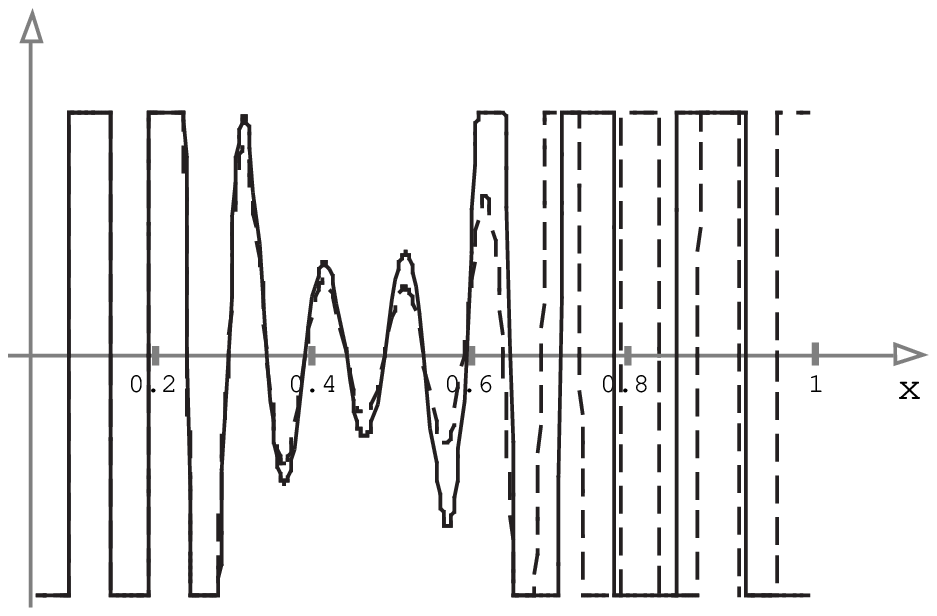}}\hskip 1cm{\epsfxsize=6cm
\epsfbox{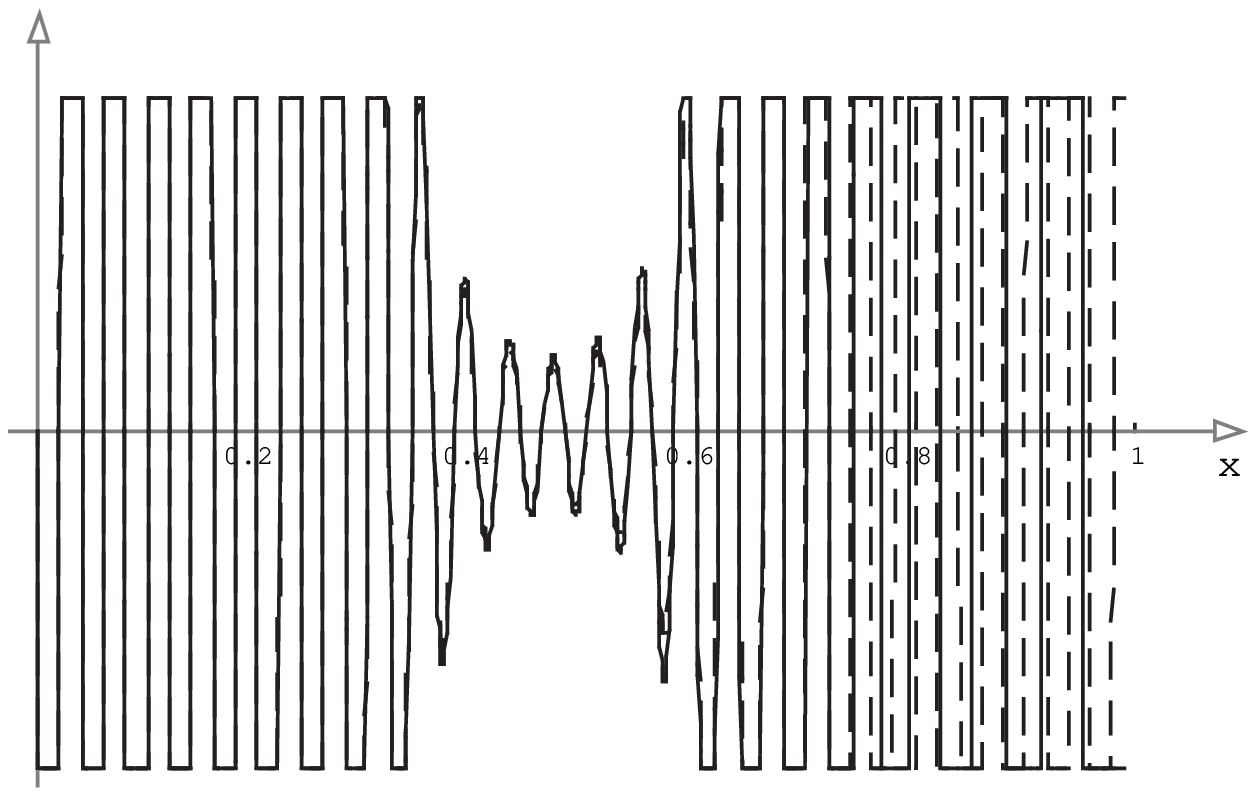}}}
\centerline{(a)\hskip 7cm (b)}
\parindent=10pt

\noindent
{\hfont Figure 3.3}. {\cfont Numerical experiment on the approximation given
in Theorem 3.1
for large $n$ and $x\in[0,1)$.
Continuous lines represent the Charlier polynomial $C_n^1(nx)$ for (a) $n=20$
and (b) $n=50$. Dashed lines represent the first order approximation given by
$e^{a/(1-x)}\Phi_0(x,n)$. Both
graphics are cut for extreme values of the polynomials.}
\vskip 2mm

\rm
\parindent=15pt

{\hfont Remark 3.4.}
When in expansion \expo\ the
expansion point $w_0$ is not equal to the saddle point $(x-1)^{-1}$,
we are not able to prove the asymptotic nature of expansion \expanchar. This
follows from the integration by parts procedure mentioned in the proof
of Theorem 3.1. On the other hand, we can show \fiasy\  directly from the
definition \ficharly\ with a change of variable like in the standard saddle
point method.

%
\subsect{1}{Details on the convergence}

It is of interest to verify the speed of convergence of the expansion
in \expanchar. We consider \ficharly, with $k=\kappa n$, and determine the
saddle point of
$(w-w_0)^\kappa (1+w)^x w^{-1}$, where $\kappa$ is large. We consider
$x$ fixed, and for $x<1$ we verify that a positive saddle point $w_+$
occurs with $w_+\sim 1/[(1-x)\kappa]$. There is a negative saddle point,
whch is not relevant. We have
$$\left[(w_+-w_0)^\kappa (1+w_+)^x w_+^{-1}\right]^n\sim
{\kappa^n e^n (1-x)^{n}\over (1-x)^{k}}.
$$
Multiplying this with $a^k n!/k!$, see \expanchar\ and \ficharly, and using in
Stirling's approximation of the factorials only the dominant parts,
that is $k!\sim (k/e)^k$, we see that
the main information on $a^k n!\Phi_k(k,n)/k!$ is given by
$$
{e^k  k^n a^k (1-x)^{n}\over (1-x)^{k} \, k^k},
$$
where $k$ is large compared with $n$ and $x$, $x < 1$. We see that the
ratio of successive terms is about $a/[(1-x)k]$.

For other values of $x$, also complex, a similar analysis can be given,
with some care in choosing the saddle points and defining the branches
in the complex plane.

%
\midinsert \hbox{\vbox{ \eightpoint {\parindent 0pt
\hfil\vbox{\offinterlineskip \hrule \halign{&\vrule#&\strut\
\hfil#\ \cr height2pt&\omit&&\omit&&\omit&&\omit&&\omit
  &&\omit&&\omit&&\omit&\cr
height2pt&\omit&&\omit&&\omit&&\omit&&\omit
  &&\omit&&\omit&&\omit&\cr
&$n$ \hfill
&&\hfill$C_n^1(0.25n)$\hfill&&${\bf C}^a_n(x,0)$\hfill
&& ${\bf C}^a_n(x,1)$\hfill&&${\bf C}^a_n(x,2)$
\hfill&&${\bf C}^a_n(x,3)$\hfill&& ${\bf C}^a_n(x,4)$
\hfill&&${\bf C}^a_n(x,5)$\hfill&\cr
height2pt&\omit&&\omit&&\omit&&\omit&&\omit
  &&\omit&&\omit&&\omit&\cr
\noalign{\hrule} height2pt&\omit&&\omit&&\omit&&\omit&&\omit
  &&\omit&&\omit&&\omit&\cr
&{10}\hfill &&-1.03630&&-0.97736&&-1.02335&&
-1.04747&&-1.00438&&-1.03633&&-1.0363&\cr
&{30}\hfill &&{4.35872}&& {4.03823}&& {
4.28867}&&{4.35077}&&{4.35762}&&{4.35858}&& {4.35870}&\cr
&{50}\hfill &&{-4.86727}&& {-4.65813}&& {
-4.82829}&&{-4.86464}&&{-4.86701}&&{-4.86725}&& {-4.86726}&\cr
&{90}\hfill &&{-2.94851}&& {-2.87926}&& {
-2.93699}&&{-2.94808}&&{-2.94848}&&{-2.94851}&& {-2.94851}&\cr
height1pt&\omit&&\omit&&\omit&&\omit&&\omit
  &&\omit&&\omit&&\omit&\cr
} \hrule}\hfil}}}
\noindent {\hfont Table 3.5.} Numerical experiment on the convergence rate of
expansion \expanchar\ for $x=.25$ and $a=1$. Here,
${\bf C}^a_n(x,N)\equiv e^{a/(1-x)}\sum_{k=0}^N{(-a)^k\over k!}\Phi_k(x,n)$,
represents the truncated series in \expanchar. All the rows are multiplied by
an appropriate constant in order to keep the numbers small.

\endinsert

\newcap{4}{Asymptotic expansions of Laguerre polynomials in terms of
Hermite polynomials}

\noindent
The Laguerre polynomials can be defined by the generating function
\eqn\genlag{
(1-t)^{-\alpha-1}e^{-tx/(1-t)}=\sum_{n=0}^\infty L_n^{\alpha}(x)t^n,
\quad \alpha,x\in\CC,\quad |t|<1.
}
and have the representation
\eqn\replag{
L_n^{\alpha}(x)=\sum_{k=0}^n(-1)^k{n+\alpha\choose n-k}\frac{x^k}{k!}.
}
For deriving the asymptotic expansion we use the Cauchy integral that
follows from \genlag:
\eqn\caulag{
L_n^{(\alpha)}(x)={1\over 2\pi i}\int_\Gamma
e^{xw/(w-1)}(1-w)^{-\alpha-1}{dw\over w^{n+1}},
}
where $\Gamma$ is a circle around the origin with radius $<1$. The
many-valued  functions $(1-w)^\mu$ appearing here and in the theorem
assume the principal branch that is equal to 1 at $w=0$

The asymptotics for large $n$, fixed $\alpha$, is considered in \frong.
For real $x$ two uniform expansions are given, one involving the $J-$Bessel
function for $x$ in an interval that contains the origin, and one in terms
of the Airy function for $x$ in an interval containing the transition
near the largest zero of $L_n^{\alpha}(x)$. In this section we give an
asymptotic expansion of $L_n^{\alpha}(nx)$ in terms of
$L_n^{1/2}(nx)$, which in fact is an Hermite polynomial. We consider $x>1$
and for these  values the expansion
is convergent  and is in particular of interest because this is the interval
that contains the large zeros and the transition point at $x=4$.

When the parameter $\alpha$ of the Laguerre
polynomial is large the asymptotic behavior can be described in terms of
Hermite polynomials (see \temlag). For example, we have the limit
\eqn\limlag{
\lim_{\alpha\to\iy}\alpha^{-n/2}L_n^\alpha(x\w{{\alpha}}+\alpha)=\frac{(-1)^n\,2^{-n/2}}{n!}
\,H_n(x/\w{{2}}).
}
In \temlopherm\ we have extended this limit by giving an asymptotic
representation for large $\alpha$ and $n$ fixed
in terms of Hermite polynomials. For more details on large $\alpha$
asymptotics we refer to \temlag. The approach in this
section is quite different because we take $\alpha$ fixed and $n$ large.

 From \caulag\   we obtain
\eqn\laguerre{
L_n^{(\alpha)}(nx)={1\over 2\pi i}\int_\Gamma f(w)
{e^{n\varphi(x,w)}\over(1-w)^{3/2}}{dw\over w},
}
where
\eqn\phifw{
\varphi(x,w)\equiv {xw\over w-1}-\log w, \quad f(w)\equiv(1-w)^{1/2  -\alpha}.
}
The function $\varphi(x,w)$ has two conjugate saddle points:
\eqn\saddles{
w^{\pm}=1-{x\over 2}\pm {i\over 2}\xi, \quad \xi=\sqrt{x(4-x)}.
}
The square root defining $\xi$ is positive for $0<x<4$; for $x\ge4$ we define
$\xi=i\sqrt{x(x-4)}$, again with positive square root.
In the expansion of the Laguerre polynomials we allow that the
saddle points coalesce.

\vfil\eject

\subsect{1}{Construction of the expansion}

The function $f(w)$ of \laguerre\
is analytic in $\Omega=\Cs\setminus[1,\infty)$ and we can expand $f(w)$ in
in a two-point Taylor expansion at the two saddle points $w^{\pm}$,
using a slightly different form of \double,

\eqn\doubleone{
f(w)=\sum_{k=0}^\infty\left\lbrack
A_k+B_kw\right\rbrack(w-w^+)^k(w-w^-)^k.
}
After substituting expansion \doubleone\ in
\laguerre\ and interchanging summation and integration we obtain
\eqn\expanlag{
L_n^{(\alpha)}(xn)=\sum_{k=0}^\infty\left\lbrack A_k\Phi_k(x,n)+B_k\Psi_k(x,n)
\right\rbrack,
}
where
\eqn\fikone{
\Phi_k(x,n)={1\over 2\pi i}\int_\Gamma(w-w^+)^k(w-w^-)^k
{e^{n\varphi(x,w)}\over(1-w)^{3/2}}{dw\over w}
}
and
\eqn\fiktwo{
\Psi_k(x,n)={1\over 2\pi i}\int_\Gamma (w-w^+)^k(w-w^-)^k
{e^{n\varphi(x,w)}\over(1-w)^{3/2}}dw.
}
We have

\eqn\fia{\eqalign{
&\Phi_0(x,n)\equiv
L_n^{(1/2)}(nx)=
{(-1)^n\over n!\,2^{2n+1}\,\sqrt{nx}}
H_{2n+1}\left(\sqrt{nx}\right),  \cr
&\Psi_0(x,n)\equiv L_{n-1}^{(1/2)}(nx)=
{(-1)^{n-1}\over (n-1)!\,2^{2n-1}\,\sqrt{nx}}
H_{2n-1}\left(\sqrt{nx}\right),\cr
}
}
and, for $k=1,2,3,...$,
\eqn\filagab{
\Phi_k(x,n)\equiv\sum_{j=0}^k\left(\matrix{k \cr j}\right)x^{k-j}
L^{(1/2-2j)}_{n-k+j}(nx), \quad
\Psi_k(x,n)\equiv\sum_{j=0}^k\left(\matrix{k \cr j}\right)x^{k-j}
L^{(1/2-2j)}_{n-k+j-1}(nx).
}
The sequences $\lbrace\Phi_k(x,n)\rbrace$ and
$\lbrace\Psi_k(x,n)\rbrace$, $k=0,1,2,\ldots$,
  satisfy the recurrences
\eqn\reculaga{
\Phi_k={1\over n-2k+3/2}\left\lbrace
a_1\Phi_{k-1}+a_2\Phi_{k-2}+
b_1\Psi_{k-1}+b_2\Psi_{k-2}
\right\rbrace,
}
$$\eqalign{
a_1&=(k-1)(x^2-2x-2)-{1\over 2}, \quad a_2=(k-1) x(2-x),\cr
b_1&=(k-1)(2-3x)+{1-x\over 2},\quad  b_2= (k-1)x(4x-x^2-2),}
$$
\eqn\reculagb{
\Psi_k=  {1\over n-2k+1/2}\left\lbrace
c_0\Phi_k+c_1\Phi_{k-1}+c_2\Phi_{k-2}+
d_1\Psi_{k-1}+d_2\Psi_{k-2},
\right\rbrace,
}
$$\eqalign{
c_0&=(2-3k)x+2(k-1)+{1-x\over 2}, \cr
c_1&=(1-k)x^3+4(k-1)x^2+kx+2(1-k)+{x-1\over2},\quad c_2=-b_2,\cr
d_1&=(4k-3)x^2+2(4-5k)x+2(k-1)+{x^2-3x+1\over 2},\cr
d_2&= (k-1)x(x^3-6x^2+9x-2).}
$$
To verify the recursions  \reculaga\  and \reculagb\ we write
\eqn\one{
\Phi_k(x,n)=-{1\over 2\pi i\,n}\int_C(w-w^+)^{k-1}(w-w^-)^{k-1}\sqrt{1-w}
{\partial e^{n\varphi(x,w)}\over\partial w}dw,
}
\eqn\fione{
\Psi_k(x,n)=-{1\over 2\pi i\,n}\int_C(w-w^+)^{k-1}(w-w^-)^{k-1}w\sqrt{1-w}
{\partial e^{n\varphi(x,w)}\over\partial w}dw,
}
\eqn\two{
\Psi_{k-1}(x,n)-w^+\Phi_{k-1}(x,n)={1\over 2\pi
i}\int_C(w-w^+)^k(w-w^-)^{k-1}{e^{n\varphi(x,w)}\over(1-w)^{3/2}} {dw\over w}
}
and
\eqn\caca{
\Psi_{k-1}(x,n)-w^-\Phi_{k-1}(x,n)={1\over 2\pi
i}\int_C(w-w^+)^{k-1}(w-w^-)^k{e^{n\varphi(x,w)}\over(1-w)^{3/2}} {dw\over w}.
}
Integrating by parts in \one\ and \fione, using \two\ and \caca\
and after straightforward manipulations we obtain \reculaga\ and \reculagb.
Formulas \fia and \filagab\ follow from \fikone\ and
\fiktwo\ after simple calculations.

\noindent
{\hfont Theorem 4.1.} {\it Expansion \expanlag\ is convergent, uniformly for
$\alpha \in \CC$ in compact sets,  and $x\ge 1+\delta>1$.
Moreover, $\{\Phi_k(x,n)\}$ and $\{\Psi_k(x,n)\}$
are asymptotic sequences for large $n$:
\eqn\lagasy{
\eqalign{
\Phi_k(x,n)=&{\cal O}\left(n^{-\lfloor (k+1)/2\rfloor}\right)\,
\left[\left|\Phi_0(x,n)\right|+\left|\Psi_0(x,n)\right|\right],\cr
\Psi_k(x,n)=&{\cal O}\left(n^{-\lfloor (k+1)/2\rfloor}\right)\,
\left[\left|\Phi_0(x,n)\right|+\left|\Psi_0(x,n)\right|\right],\cr
}
}
as $n\to\infty$ and $k=0,1,2,...$.
}

\noindent
{\it Proof.\quad}
We apply  Proposition 2.2. Expansion \doubleone\ is uniformly convergent
for $w$ inside the Cassini oval with foci $w^+$ and $w^-$ and
"radius" $r=\vert w_1-w^+\vert\vert w_1-w^-\vert$, where $w_1=1$ is the
singular point of $f$. Using \saddles\ it follows that $r=x$.
The points $w$ are inside
the Cassini oval if they satisfy $\vert w-w^+\vert\vert w-w^-\vert<r=x$.
Because $w^+w^-=1$, the origin $w=0$ is
inside the oval only if $x>1$. Hence, the contour $\Gamma$ of
\laguerre\ can be taken  completely inside the oval only if $x>1$ (see
also Figure 4.2).
This proves the convergence of \expanlag\  for $x>1$.
The asymptotic behavior in \lagasy\  follows from \fia\ and the recursions
\reculaga\ and \reculagb. More detailed asymptotic information can be
obtained from  the integrals in \fikone\ and \fiktwo.
\hfill $\boxe$

\bigskip

\centerline{{\epsfxsize=9cm \epsfbox{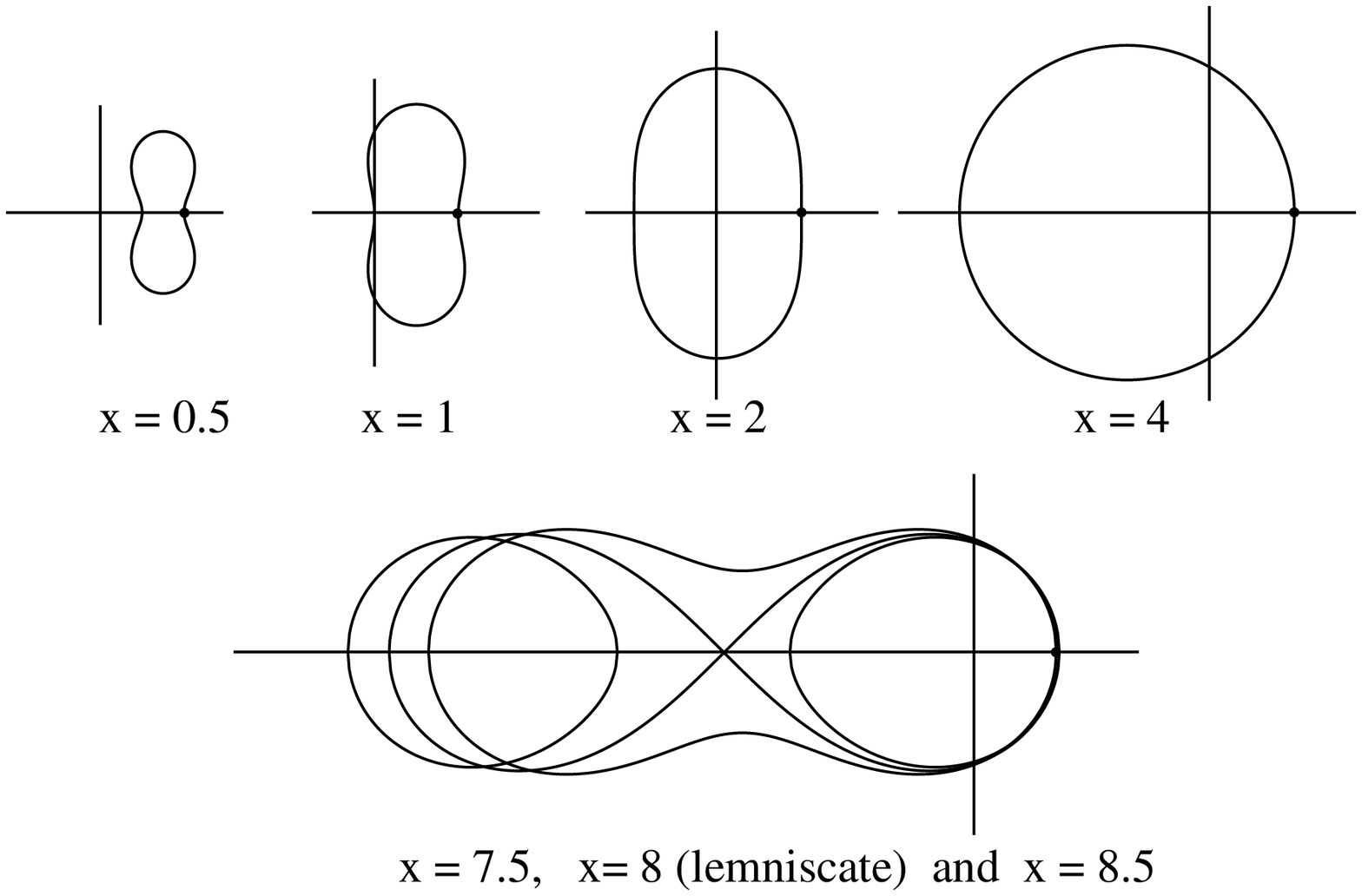}}}
\parindent=10pt

\noindent
{\hfont Figure 4.2.} {\cfont Cassini ovals for the expansion \doubleone\ for
several values of $x$. For $w$ inside the ovals, the expansion is
convergent. For $0<x<1$ the origin is outside the oval; for $x=4$
it is a circle, for $x=8$ a lemniscate. For $x>8$ the oval splits up into
two parts. All ovals go through the point $w=1$, a singular point of $f(w)$.}

\rm
\parindent=15pt


%
\centerline{{\epsfxsize=6cm \epsfbox{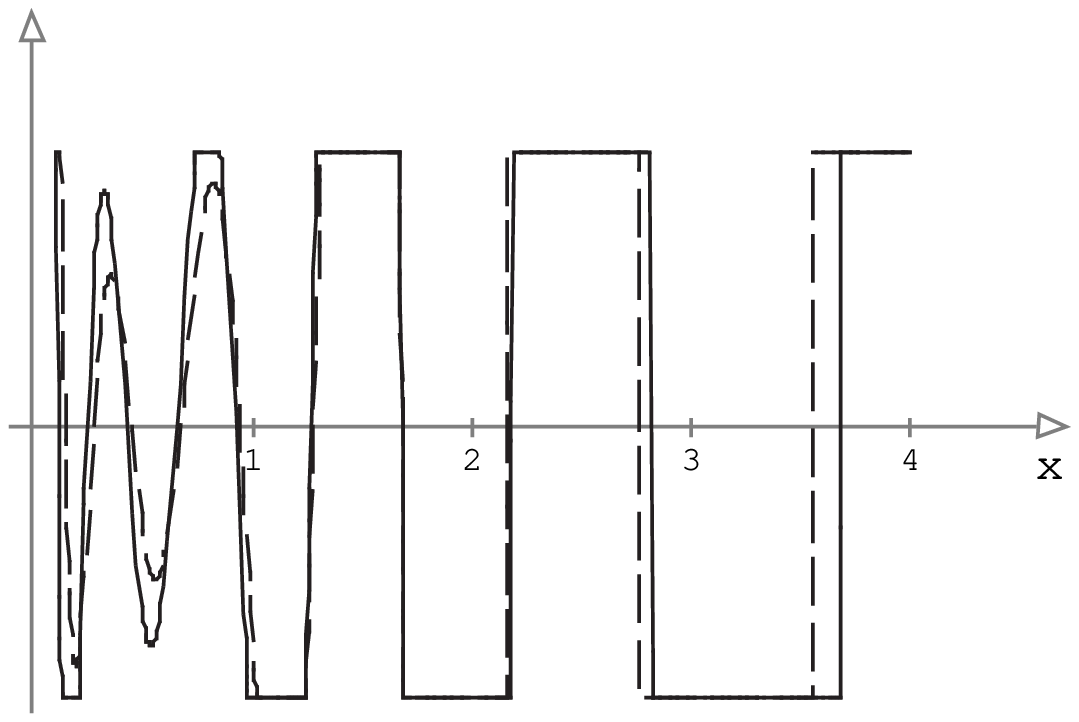}}\hskip 2cm
{\epsfxsize=6cm \epsfbox{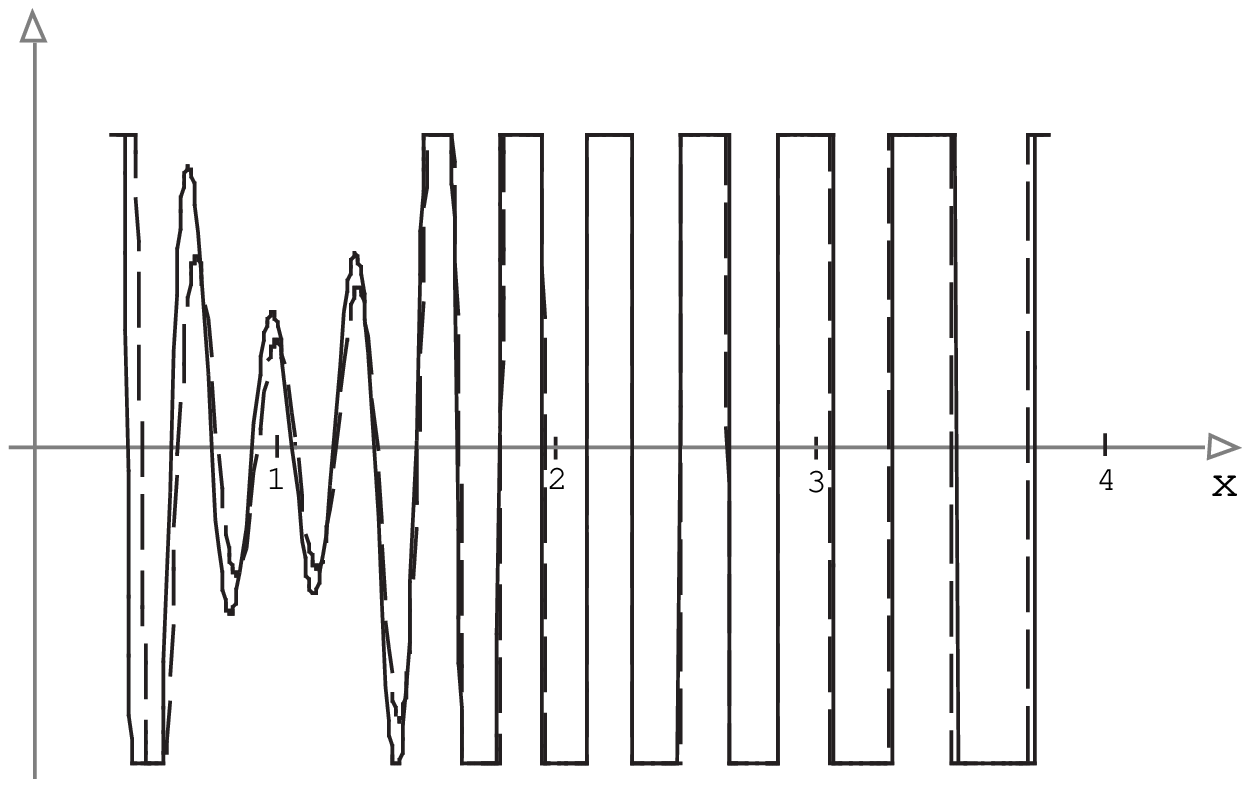}}}
\centerline{(a) \hskip 7cm (b)}
\parindent=10pt

\noindent
{\hfont Figure 4.3.} {\cfont  Numerical experiments on the approximation of
Theorem 4.2 for large $n$ and $x\in(0,4]$.
Continuous lines represent the Laguerre polynomial $L_n^{(4)}(nx)$ for (a)
$n=10$
and (b) $n=20$. Dashed lines represent the first order approximation given by
$A_0\Phi_0(x,n)+B_0\Psi_0(x,n)$. Both
graphics are cut for extreme values of the polynomials.}

\rm
\parindent=15pt

\noindent
{\hfont Remark 4.4.} The expansion in \expanlag\ has a meaning for all complex
$x$, and has for all fixed $x$ an asymptotic meaning. The expansion is
uniformly convergent for $\vert x\vert\ge 1+\delta>1$.
%


\subsect{2}{Details on the coefficients}

The expressions \coefibis\ and \aprime\ for the coefficients in
\double\ can be used in the present case also. We first write
\eqn\doublew{
f(w)=\sum_{k=0}^\infty\left\lbrack a_k(w-w^-)+a_k'(w-w^+)\right\rbrack
[(w-w^-)(w-w^+)]^k,
}
and compare this with \doubleone. By comparing coefficients of equal
powers, it follows that
$A_k$ and $B_k$ can be expressed in terms of $a_k$ and $a'_k$. We have
for $k=0,1,2,\ldots$
$$
A_k=-a_kw^+-a_k'w^-,\quad  B_k=a_k+a_k'.
$$
We have
\eqn\abk{
A_k\equiv -2\Re\left(w^+
a_k\right),\quad
B_k\equiv 2\Re\left(a_k\right),\quad
a_0=i(1-w^-)^{1/2-\alpha} \xi^{-1},
}
and, for $k=1,2,3,...$,
\eqn\ak{
a_k= \sum_{j=0}^k{(k+j-1)!(\alpha-1/2)_{k-j}\over
  k!j!(k-j)!(i\xi)^{k+j+1}}
\left\{ {(-1)^jj\over
(1-w^+)^{\alpha+k-j-1/2}}-
{(-1)^kk\over (1-w^-)^{\alpha+k-j-1/2}}
\right\}.
}

The coefficients can also be computed from the recursion relations
\eqn\recakbk{\eqalign{
&x(k+1)A_{k+1}-x(k+1)B_{k+1}       =(\alpha-1/2+2k)A_k-(xk+1)B_k,\cr
&x(k+1)A_{k+1}-x(x-3)(k+1)B_{k+1}  =(\alpha+1/2+2k)B_k,\cr
  }}
where $k=0,1,2,\ldots\ $. Let, for $1<x\le4$, $x=4\sin^2(\theta/2)$. Then
$w^{\pm}=e^{\pm i\theta}$, and
$$
A_0=-2^{\beta+1} \sin^\beta(\theta/2)\cos[(\theta-\pi)\beta/2-\theta],\quad
B_0= 2^{\beta+1} \sin^\beta(\theta/2)\cos[(\theta-\pi)\beta/2].
$$
where $\beta=1/2-\alpha$. This gives real expressions for the first
coefficients to start the recursion relations in \recakbk. For $x\ge4$ we
can obtain expressions in terms of hyperbolic functions by writing
$x=4\cosh^2(\theta/2)$, which gives $w^{\pm}=-e^{\pm\theta}$ and
$$
A_0=-2^{\beta+1} \cosh^\beta(\theta/2)\cosh[(\theta(\beta/2-1)],\quad
B_0= 2^{\beta+1} \cosh^\beta(\theta/2)\cosh(\theta\beta/2).
$$


\subsect{3}{An alternative form of the expansion}

By using in \laguerre\  the substitution
\eqn\fsubs{
f(w)=\alpha_0 +\beta_0 w +(w-w^-)(w-w^+) g_0(w),
}
where $\alpha_0$ and $\beta_0$ follow from substituting $w=w^{\pm}$,
we obtain by integrating by parts
\eqn\laltf{
L_n^{(\alpha)}(nx)=\alpha_0\Phi_0(x,n)+\beta_0\Psi_0(x,n)+
{1\over 2\pi i\,n} \int_{\Gamma}
f_1(w)
{e^{n\varphi(x,w)}\over(1-w)^{3/2}}{dw\over w},
}
where
\eqn\fone{
f_1(w)=w(1-w)^{3/2}{d\over dw} \left[(\sqrt{1-w}\,g_0(w)\right].
}
Continuing this procedure we obtain the expansion in negative powers of
$n$:
\eqn\lalt{
L_n^{(\alpha)}(nx)=\Phi_0(x,n)\sum_{k=0}^\infty {\alpha_k\over n^k}
+\Psi_0(x,n)\sum_{k=0}^\infty {\beta_k\over n^k},
}
where $\alpha_k$ and $\beta_k$ follow from
\eqn\abk{
f_k(w^-)=\alpha_k+ \beta_kw^-,\quad f_k(w^+)=\alpha_k+ \beta_kw^+,
}
where with $f_0(w)=f(w)$ and for $k=0,1,2,\ldots$
\eqn\fk{
\eqalign{
f_k(w)&=\alpha_k+ \beta_kw+(w-w^-)(w-w^+) g_k(w),\cr
f_{k+1}(w)&=w(1-w)^{3/2}{d\over dw} \left[(\sqrt{1-w}\,g_k(w)\right].\cr
}
}
The expansion in \lalt\ also follows from re-arranging expansion
\expanlag\ by using the recursion relations for
$\Phi_k(x,n)$ and $\Psi_k(x,n)$ in \reculaga\  and \reculagb.
%
%
%
\midinsert \hbox{\vbox{ \eightpoint {\parindent 0pt
\hfil\vbox{\offinterlineskip \hrule \halign{&\vrule#&\strut\
\hfil#\ \cr height2pt&\omit&&\omit&&\omit&&\omit&&\omit
  &&\omit&&\omit&&\omit&\cr
height2pt&\omit&&\omit&&\omit&&\omit&&\omit
  &&\omit&&\omit&&\omit&\cr
&$n$ \hfill
&&\hfill$L_n^{(1)}(3.5n)$\hfill&&${\bf L}^{(1)}_n(3.5,0)$\hfill
&& ${\bf L}^{(1)}_n(3.5,1)$\hfill&&${\bf L}^{(1)}_n(3.5,2)$
\hfill&&${\bf L}^{(1)}_n(3.5,3)$\hfill&& ${\bf L}^{(1)}_n(3.5,4)$
\hfill&&${\bf L}^{(1)}_n(3.5,5)$\hfill&\cr
height2pt&\omit&&\omit&&\omit&&\omit&&\omit
  &&\omit&&\omit&&\omit&\cr
\noalign{\hrule} height2pt&\omit&&\omit&&\omit&&\omit&&\omit
  &&\omit&&\omit&&\omit&\cr
&{10}\hfill &&0.340506&&0.343249&&0.341724&&
0.340495&&0.340449&&0.340490&&0.340504&\cr
&{30}\hfill &&-8.94039&&-8.86531&&-9.03530
&&-8.95798&&-8.94213&&-8.94045&&-8.94038&\cr
&{50}\hfill &&-5.05678&&-5.05941&&-5.06764
&&-5.05801&&-5.05689&&-5.05680&&-5.05678&\cr
&{90}\hfill &&6.56556&&6.56328&&6.57572
&&6.56601&&6.56547&&6.56553&&6.56556&\cr
height1pt&\omit&&\omit&&\omit&&\omit&&\omit
  &&\omit&&\omit&&\omit&\cr
} \hrule}\hfil}}}
\noindent {\hfont Table 4.5.}\ Numerical experiment on the convergence rate of
expansion \expanlag\ for $x=3.5$ and $a=1$. Here,
${\bf L}^{(\alpha)}_n(x,N)\equiv
\sum_{k=0}^N[A_k\Phi_k(x,n)+B_k\Psi_k(x,n)]$,
represents the truncated series in \expanlag. All the rows are multiplied by
an appropriate constant in order to keep the numbers small.

\endinsert

\newcap{5}{Asymptotic expansions of Jacobi polynomials in terms of Chebyshev
polynomials}

\noindent
The large $n$ asymptotics for the Jacobi polynomials is discussed in
[\bateman, Vol. II, \S 10.14], in particular for $x\in (-1,1)$.
For $x$ bounded away from the points $\pm 1$ elementary functions
(sine and cosine functions) can be used for describing the asymptotics.
For $x$ close to $\pm1$ Bessel functions can be used (Hilb-type formulas).

In this section we develop a convergent expansion that is valid for
$x\in(-1,1)$ and the terms of the expansion constitute asymptotic
scales for large $n$. It is possible to extend the results to complex
values of $x$, but this will not be considered here.
The first approximants are Chebyshev polynomials,
which in fact are elementary functions, and the other terms
can be obtained from recursions that show the asymptotic property.

\vfil\eject

\subsect{1}{Construction of the expansion}

Starting point is the integral representation
that follows from  [\bateman, Vol. II, p. 172]:
\eqn\jacobi{
P_n^{(\alpha,\beta)}(x)={1\over 2\pi i}{(-1)^n\over 2^n}
\int_\Gamma
{(1-w-x)^{\alpha}(1+w+x)^{\beta}\over
(1-x)^\alpha(1+x)^\beta}e^{n\varphi(x,w)}{dw\over w},
}
where we consider $x\in(-1,1)$; the function $\varphi(x,w)$ is defined by
$$
\varphi(x,w)\equiv \log(1+w+x)+\log(1-w-x)-\log w
$$
and $\Gamma$ is
a simple closed contour, in the positive sense, around $w=0$. The 
points $w=-x\pm 1$
are outside the contour, and  
$(1-w-x)^{\alpha}/(1-x)^\alpha$ and $(1+w+x)^{\beta}/
(1+x)^\beta$ are to be taken as unity when $w=0$.

The function $\varphi(x,w)$ has two conjugate saddle points:
\eqn\sadd{
w^{\pm}=\pm w_0, \quad w_0=i\sqrt{1-x^2}.
}
We expand the integral by using the function
$$
f(w)\equiv(1-x)^{-\alpha-{1\over 2}}(1+x)^{-\beta-{1\over 2}}
(1-w-x)^{\alpha+{1\over 2}}(1+w+x)^{\beta+{1\over 2}},
$$
which is analytic in
$$
\Cs\setminus\lbrace(-\infty,-1-x]\bigcup[1-x,\infty)\rbrace.
$$
We expand, using a slightly different form of \double,
\eqn\doubletwo{
f(w)=\sum_{k=0}^\infty
\left\lbrack A_k+B_k\,w\right\rbrack(w^2-w_0^2)^k
}
where the coefficients $A_k$ and $B_k$ can be expressed in
terms of the derivatives of
$f(w)$ at $w=\pm w_0$; see the next subsection for more details.

After substituting expansion \doubletwo\ in
\jacobi\ and interchanging summation and integration we obtain
\eqn\expanjac{
P_n^{(\alpha,\beta)}(x)=\sum_{k=0}^\infty\left\lbrack
A_k\Phi_k(x,n)+B_k\Psi_k(x,n)
\right\rbrack,
}
where
\eqn\phik{
\Phi_k(x,n)={(-1)^n\over 2\pi i}{\sqrt{1-x^2}\over 2^n}
\int_\Gamma{(w^2+1-x^2)^k\over
W(x,w)}e^{n\varphi(x,w)}{dw\over w},
}
\eqn\psik{
\Psi_k(x,n)={(-1)^n\over 2\pi i}{1\over 2^n\sqrt{1-x^2}}
\int_\Gamma(w^2+1-x^2)^k
W(x,w)\,e^{n\varphi(x,w)}{dw\over w},
}
and
\eqn\wxw{
W(x,w)\equiv\sqrt{(1-w-x)(1+w+x)}.
}
We have
\eqn\fipsizero{
\Phi_0(x,n)\equiv P_n^{(-1/2,-1/2)}(x), \quad
\Psi_0(x,n)\equiv -{1\over 2}(1-x^2)\,P_{n-1}^{(1/2,1/2)}(x).
}
These Jacobi polynomials are Chebyshev
polynomials:
\eqn\psizero{
P_n^{(-1/2,-1/2)}(x)
={2^{-2n}(2n)!\over (n!)^2} T_n(x),\quad
P_{n}^{(1/2,1/2)}(x)={2^{-2n}(2n+1)!\over n!\,(n+1)!} U_n(x).
}
In terms of elementary functions:
\eqn\cheb{
T_n(\cos\theta)=\cos n\theta,\quad
U_n(\cos\theta)={\sin(n+1)\theta \over\sin\theta}.
}
Furthermore
\eqn\fiunocero{
\Phi_1(x,n)\equiv
{1-x^2\over4(n+1)}\left[4(n+1)\Phi_0(x,n)-
4x(n-1)\Psi_0(x,n)+(2n-1)\Psi_0(x,n-1)\right],
}
\eqn\fiunouno{
\Psi_1(x,n)\equiv{-3(1-x^2)\over 4(n+1)(n+2)}\left\{
2[n+1-2x^2n]\Psi_0(x,n)+x(2n-1)\Psi_0(x,n-1)\right\},
}
and, for $k=2,3,4,...$,
\eqn\fikjaca{
\Phi_k(x,n)\equiv\sum_{j=0}^k\left(\matrix{k\cr
j}\right){(1-x^2)^{k+j}\over 4^j}
P_{n-2j}^{(2j-1/2,2j-1/2)}(x),
}
\eqn\fikcacb{
\Psi_k(x,n)\equiv -{1\over 2}(1-x^2)\sum_{j=0}^k\left(\matrix{k\cr
j}\right){(1-x^2)^{k+j}\over 4^j}
P_{n-1-2j}^{(2j+1/2,2j+1/2)}(x).
}
The sequences $\lbrace\Phi_k(x,n)\rbrace$ and
$\lbrace\Psi_k(x,n)\rbrace$
satisfy the recursion relations
\eqn\recujaca{
\Phi_k=  {1\over n+2k-1}\left\lbrack a_1\Phi_{k-1}+a_2\Phi_{k-2}
+b_1\Psi_{k-1} +b_2\Psi_{k-2}\right\rbrack,
}
\eqn\recujacb{
\Psi_k=  {1\over n+2k}\left\lbrack c_0 \Phi_k + c_1\Phi_{k-1}+c_2\Phi_{k-2}
+d_1\Psi_{k-1} +d_2\Psi_{k-2}\right\rbrack,
}
where
$$
\eqalign{
a_1&=(1-x^2)(6k-5), \quad a_2=-4(k-1)(1-x^2)^2,\cr
b_1&=-x(4k-3),\quad b_2= 4 x (k-1) (1-x^2),\cr
c_0&=-x(4k-1),\quad c_1= x(1-x^2)(8k-5),\quad c_2=-4x(1-x^2)^2(k-1),\cr
d_1&=3(2k-1)(1-x^2),\quad d_2=-4 (1-x^2)^2 (k-1).\cr
}
$$
For the relations between Jacobi and Chebyshev polynomials we refer to [\temsf\
p. 152-153]. The expressions in \fiunocero\ and \fiunouno\ follow from
contiguous relations  of the Jacobi polynomials (see[\temsf\ p. 166]).

To verify the recursions in \recujaca\ and \recujacb, we write
\eqn\ones{
\Phi_k(x,n)={-1\over 2\pi i}{(-1)^n\over n\,2^n\sqrt{1-x^2}}
\int_\Gamma(w^2+1-x^2)^{k-1}\,W(x,w)\,
{\partial e^{n\varphi(x,w)}\over\partial w}\,dw,
}
\eqn\twos{
\Psi_k(x,n)={-1\over 2\pi i}{(-1)^n\over n\,2^n\sqrt{1-x^2}}
\int_\Gamma w(w^2+1-x^2)^{k-1}\,W(x,w)\,
{\partial e^{n\varphi(x,w)}\over\partial w}\,dw,
}
Integrating by parts in \ones\ and \twos\ and
after straightforward manipulations we obtain \recujaca\ and \recujacb.
The asymptotic behavior pointed out above follows from
the definition of the Jacobi polynomials and the recurrences \recujaca\ and
\recujacb.

\noindent
{\hfont Theorem 5.1.} {\it Expansion \expanjac\ is convergent for $x\in(-1,1)$.
Moreover, $\{\Phi_k(x,n)\}$ and $\{\Psi_k(x,n)\}$
are asymptotic sequences for large $n$:
$$\Phi_k(x,n)={\cal O}\left(n^{-\lfloor (k+1)/2\rfloor}\right),\quad
\Psi_k(x,n)={\cal O}\left(n^{-\lfloor (k+1)/2\rfloor}\right),$$
as $n\to\infty$ and $k=0,1,2,...$.
}

\noindent
{\it Proof.\quad}
Expansion \doubletwo\
is uniformly convergent for $w$ inside the Cassini oval of focus $w^\pm$ and
"radius"
$$
r=\min \lbrace\vert 1+x+w^+\vert\vert 1+x+w^-\vert,
\vert 1-x+w^+\vert\vert 1-x+w^-\vert\rbrace=
2(1-\vert x\vert).
$$
Hence, \doubletwo\ is convergent for $w\in\Cs$ such that
$\vert w^2+1-x^2\vert<2(1-\vert x\vert)$.
 From Proposition 2.2 it follows that expansion \expanjac\ is convergent for
$1-x^2=\vert w^+w^-\vert<2(1-\vert x\vert)$, that is, $\forall$
$x\in(-1,1)$.
The asymptotic property follows from the recursion relations in
\recujaca\ and \recujacb.
\hfill$\boxe$

\subsect{2}{Details on the coefficients}

The expressions for  \coefibis\ and \aprime\ for the coefficients in
\double\ can be used in the present case also. We first write
\eqn\doublew{
f(w)=\sum_{k=0}^\infty\left\lbrack a_k(w-w_0)+a_k'(w+w_0)\right\rbrack
(w^2-w_0^2)^k,
}
and compare this with \doubletwo. By comparing coefficients of equal
powers, it follows that
$A_k$ and $B_k$ can be expressed in terms of $a_k$ and $a'_k$. We have
for $k=0,1,2,\ldots$
$$
A_k=(a_k'-a_k)w_0,\quad  B_k=a_k+a_k'.
$$
After straightforward manipulations we find
\eqn\abk{
A_k\equiv {2\Im\left\lbrack a_k\sqrt{1-x^2}\right\rbrack\over
(1-x)^{\alpha+1/2}(1+x)^{\beta+1/2}}
\hskip 1cm {\rm and} \hskip 1cm
B_k\equiv {2\Re a_k \over (1-x)^{\alpha+1/2}(1+x)^{\beta+1/2}},
}
where
%
\eqn\ajo{
a_0=-{\left(1-x+w_0\right)^{\alpha+1/2}\left(1+x-w_0\right
)^{\beta+1/2}
\over 2w_0}
}
and, for $m=1,2,3,...$,
\eqn\ajk{
\eqalign{
a_m=&{1\over m!}\sum_{k=0}^m{(m+k-1)!\over
k!(m-k)!(2w_0)^{m+k+1}}\times
\cr &
\sum_{j=0}^{m-k}\left(\matrix{m-k \cr j \cr}\right)\left\lbrace
{k(-1)^{m+j}(-\alpha-1/2)_j(-\beta-1/2)_{m-k-j}\over
\left(1-x-w_0\right)^{j-\alpha-1/2}
\left(1+x+w_0\right)^{m-k-j-\beta-1/2}}- \right.\cr &\left.
{m(-1)^{k+j}(-\alpha-1/2)_j(-\beta-1/2)_{m-k-j}\over
\left(1-x+w_0\right)^{j-\alpha-1/2}
\left(1+x-w_0\right)^{m-k-j-\beta-1/2}}
\right\rbrace.\cr}
}
\bigskip
\centerline{{\epsfxsize=6cm \epsfbox{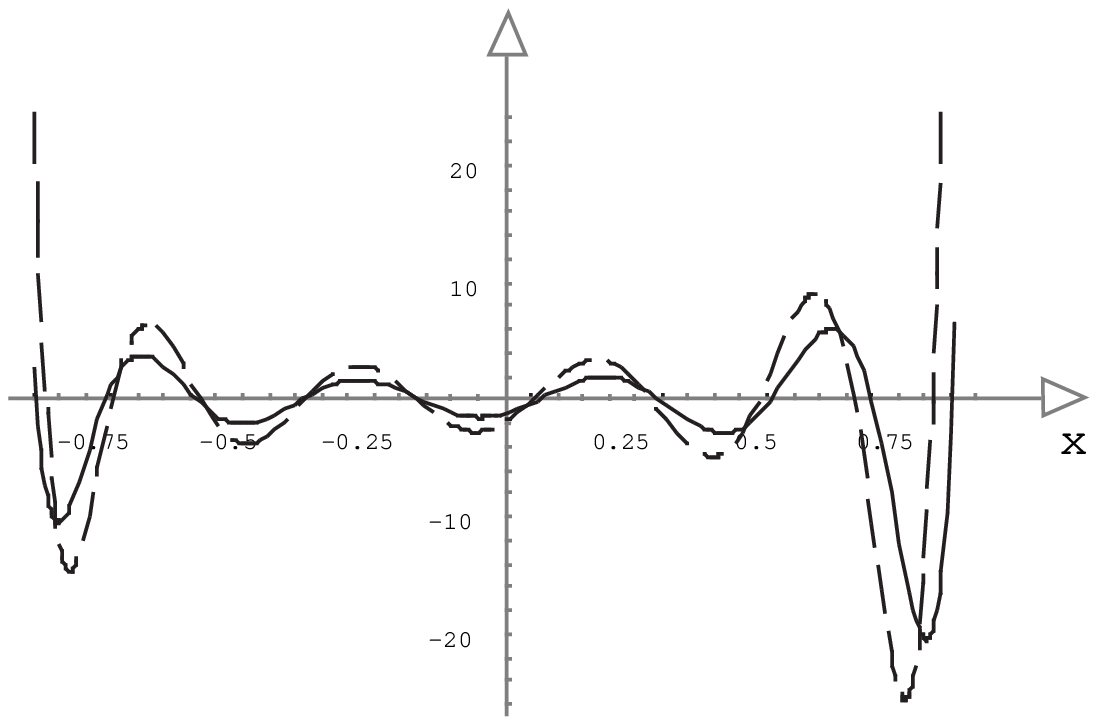}}
\hskip 1cm{\epsfxsize=6cm \epsfbox{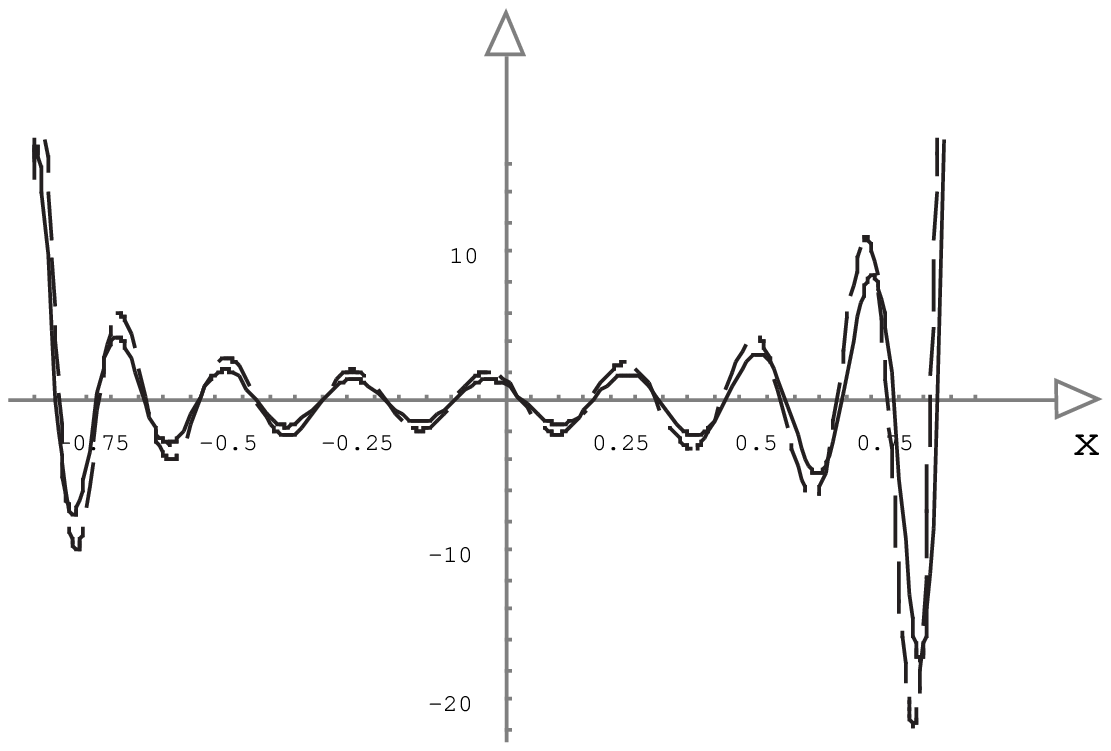}}}
\centerline{(a) \hskip 7cm (b)}
\parindent=10pt

\noindent
{\hfont Figure 5.2.} {\cfont Numerical experiments about the approximation
of Theorem 5.1
for large $n$ and $-1<x<1$.
Continuous lines represent the Jacobi polynomial $P_n^{(3,4)}(x)$ for 
(a) $n=10$
and (b) $n=20$. Dashed lines represent the first order approximation given by
$A_0\Phi_0+B_0\Psi_0$.}
\vskip 2mm

\rm
\parindent=15pt

%
\midinsert \hbox{\vbox{ \eightpoint {\parindent 0pt
\hfil\vbox{\offinterlineskip \hrule \halign{&\vrule#&\strut\
\hfil#\ \cr height2pt&\omit&&\omit&&\omit&&\omit&&\omit
  &&\omit&&\omit&&\omit&\cr
height2pt&\omit&&\omit&&\omit&&\omit&&\omit
  &&\omit&&\omit&&\omit&\cr
&$n$ \hfill
&&\hfill$P_n^{(3/2,1/2)}(0)$\hfill&&${\bf P}_n(0,0)$\hfill
&& ${\bf P}_n(0,1)$\hfill&&${\bf P}_n(0,2)$
\hfill&&${\bf P}_n(0,3)$\hfill&& ${\bf P}_n(0,4)$
\hfill&&${\bf P}_n(0,5)$\hfill&\cr
height2pt&\omit&&\omit&&\omit&&\omit&&\omit
  &&\omit&&\omit&&\omit&\cr
\noalign{\hrule} height2pt&\omit&&\omit&&\omit&&\omit&&\omit
  &&\omit&&\omit&&\omit&\cr
&{10}\hfill &&-0.336376&&-0.348029&&-0.348029
&&-0.337153&&-0.336376&&-0.336340&&-0.336360&\cr
&{30}\hfill &&-0.201847&&-0.204304&&-0.204304
&&-0.201909&&-0.201839&&-0.201845&&-0.201847&\cr
&{50}\hfill &&-0.157618&&-0.158781&&-0.158781
&&-0.157636&&-0.157615&&-0.157617&&-0.157618&\cr
&{90}\hfill &&-0.118124&&-0.118612&&-0.118612
&&-0.118128&&-0.118123&&-0.118124&&-0.118124&\cr
height1pt&\omit&&\omit&&\omit&&\omit&&\omit
  &&\omit&&\omit&&\omit&\cr
} \hrule}\hfil}}}
\noindent {\hfont Table 5.3.} Numerical experiment on the convergence rate of
expansion \expanjac\ for $x=0$, $\alpha=3/2$ and $\beta=1/2$. Here,
${\bf P}_n(x,N)\equiv
\sum_{k=0}^N[A_k\Phi_k(x,n)+B_k\Psi_k(x,n)]$,
represents the truncated series in \expanjac.

\endinsert

\vfill\eject

\newcap{6}{A few other examples of convergent asymptotic expansions}

\noindent
Consider the integral for the modified Bessel function
(for properties of the special functions in this section we refer to
\bateman\  or \temsf)
\eqn\ibz{
I_\nu(z)={(2z)^\nu e^z\over\sqrt{\pi}\Gamma(\nu+{1\over2})}\,
\int_0^1 e^{-2zt} [t(1-t)]^{\nu-{1\over 2}}\,dt.
}

The usual method for obtaining the asymptotic expansion for large $z$
consists of substituting the expansion $(1-t)^{\nu-{1\over2}}=\sum
{\nu-{1\over2}\choose k}  (-t)^k$
and  by interchanging the order of summation and integration.
The resulting integrals are not evaluated over $[0,1]$ but
over $[0,\infty)$. This latter step gives the divergent asymptotic
expansions. If we integrate over $[0,1]$ we obtain an expansion
in terms of incomplete gamma functions:
\eqn\ibex{
I_\nu(z)={e^z\over\sqrt{2z\pi}\Gamma(\nu+{1\over2})}\,
\sum_{k=0}^\infty (-1)^{k}{\nu-{1\over2}\choose k}
{\gamma(\nu+k+1/2,2z)\over(2z)^k}.
}
For fixed values of $z$ the incomplete gamma functions
have the asymptotic behaviour
\eqn\igas{
\gamma(\nu+k+1/2,2z)={e^{-2z}(2z)^{\nu+k+\frac12}\over \nu+k+\frac12}\,
\left[1+{\cal O}(k^{-1})\right],\quad k\to\infty,
}
and we see that the terms behave like ${\cal O}(k^{-\nu-3/2})$. So,
convergence is guaranteed if $\Re\nu>-{1\over 2}$.
A further examination of the terms of the expansion shows that
for large $z$ it is better to use the asymptotic property of
the expansion (the ratio of successive terms is of order ${\cal O}(1/z)$).
The incomplete gamma functions can be computed by using
a backward recusion scheme.

Another example  is the $K-$Bessel function given by
\eqn\kbz{
K_\nu(z)={\sqrt{\pi}(2z)^\nu e^{-z}\over\Gamma(\nu+{1\over2})}\,
\int_0^\infty e^{-2zt} [t(1+t)]^{\nu-{1\over 2}}\,dt.
}
Again expanding $(t+1)^{\nu-1/2}=\sum{\nu-{1\over2}\choose k}  (t)^k$
gives the standard expansion. A convergent expansion can be obtained by
using
\eqn\kbz{
(t+1)^{\nu-1/2}=\sum c_k  \left({t\over t+1}\right)^k, \quad
c_k=(-1)^k{{1\over2}-\nu\choose k}.
}
This gives a convergent expansion in terms of confluent 
hypergeometric functions
\eqn\kbex{
K_\nu(z)=\sqrt{\pi\over 2z}e^{-z}
\sum_{k=0}^\infty (-1)^{k}{(\nu+{1\over 2})_k (\nu-{1\over 2})_k\over k!}
U\left(k,\tfrac12-\nu,2z\right).
}
We have
\eqn\ufun{
U\left(0,\tfrac12-\nu,2z\right)=1,\quad
U\left(1,\tfrac12-\nu,2z\right)=(2z)^{\nu+{1\over2}}\,e^{2z}\,
\Gamma\left(-\tfrac12-\nu,2z\right),
}
again an incomplete gamma function. Other $U-$functions can be obtained by
recursion. For large $k$ the $U-$function behaves like \luke
$$
k!\,U\left(k,\tfrac12-\nu,2z\right)=
{\cal O}\left(k^\alpha\,e^{-2\sqrt{2kz}}\right),
$$
where $\alpha$ is some constant.
It follows that the convergence is better than in the previous example.

As a final example we consider an expansion of the Kummer function.
Tricomi \tric\ has derived several convergent expansions of the
${}_1F_1-$function in terms of Bessel functions that are useful
for evaluating the function when the parameters are large. For example, we have
\tric
\eqn\trico{
{}_1F_1(a,c;z)=e^{{1\over2}z}\G(c)(\kappa z)^{(1-c)/2}\sum_{n=0}^\infty\,
A_n(\kappa,\tfrac12c)\,\left({z\over4\kappa}\right)^{n/2}\,
J_{c-1+n}\left(2\sqrt{{\kappa z}}\right),
}
where $\kappa=\frac12c-a$ and the $A_n(\kappa,\lambda)$
are coefficients in the generating function
\eqn\tricoef{e^{2\kappa z}\,(1-z)^{\kappa-\lambda}\,(1+z)^{-\kappa-\lambda}=
\sum_{n=0}^\infty\,A_n(\kappa,\lambda)\,z^n.
}
The series in \trico\  is convergent in the entire $z-$plane. Moreover, it can
be used for the evaluation of ${}_1F_1(a,c;z)$ for large $\kappa$,
because the series has an asymptotic property. For further details
on these expansions we refer to \tric.

\newcap{7}{Acknowledgements}

\noindent
The authors thank the referee for the comments on the 
first version of the paper.

J. L. L\'opez wants to thank the C.W.I. of Amsterdam for its scientific and
financial support
during the realization of this work. The financial support of the saving
bank {\it Caja Rural de
Navarra} is also acknowledged.

{\ninepoint
\listrefs
}
\end

%% file: 102153macro.tex
\newdimen\fullhsize \newdimen\hstitle \newdimen\hsbody
\tolerance=1000\hfuzz=2pt
%
\magnification=\magstephalf
\baselineskip=14pt plus 2pt minus 1pt
\parskip=4pt
\parindent=.31truein
\hoffset=.6truein
\voffset=.8truein
\hsbody=\hsize \hstitle=\hsize 
\hsize=12.6truecm
\vsize=19.8truecm
\catcode`\@=11 
\newcount\yearltd\yearltd=\year\advance\yearltd by -1900

%
%
\def\draftmode{\def\draftdate{{\rm preliminary draft:
\number\month/\number\day/\number\yearltd\ \ \hourmin}}%
\headline={\hfil\draftdate}\writelabels\baselineskip=14pt plus 2pt minus 2pt
{\count255=\time\divide\count255 by 60 \xdef\hourmin{\number\count255}
	\multiply\count255 by-60\advance\count255 by\time
    \xdef\hourmin{\hourmin:\ifnum\count255<10 0\fi\the\count255}}}

\def\draftdate{}
\def\nolabels{\def\eqnlabel##1{}\def\eqlabel##1{}\def\reflabel##1{}}
\def\writelabels{\def\eqnlabel##1{%
{\escapechar=` \hfill\rlap{\hskip.09in\string##1}}}%
\def\eqlabel##1{{\escapechar=` \rlap{\hskip.09in\string##1}}}%
\def\reflabel##1{\noexpand\llap{\string\string\string##1\hskip.31in}}}
\nolabels
%
\global\newcount\secno \global\secno=0
\global\newcount\meqno \global\meqno=1
\def\appendix#1#2{\global\meqno=1\xdef\secsym{\hbox{#1.}}\bigbreak\bigskip
\noindent{\bf Appendix #1. #2}\par\nobreak\medskip\nobreak}
%
%
\def\eqnn#1{\xdef #1{(\secsym\the\meqno)}%
\global\advance\meqno by1\eqnlabel#1}
\def\eqna#1{\xdef #1##1{\hbox{$(\secsym\the\meqno##1)$}}%
\global\advance\meqno by1\eqnlabel{#1$\{\}$}}
\def\eqn#1#2{\xdef #1{(\secsym\the\meqno)}\global\advance\meqno by1%
$$#2\eqno#1\eqlabel#1$$}
%
%
\global\newcount\refno \global\refno=1
\newwrite\rfile
\def\ref{[\the\refno]\nref}
\def\nref#1{\xdef#1{[\the\refno]}\ifnum\refno=1\immediate
\openout\rfile=refs.tmp\fi\global\advance\refno by1\chardef\wfile=\rfile
\immediate\write\rfile{\noexpand\item{#1\ }\reflabel{#1}\pctsign}\findarg}
\def\findarg#1#{\begingroup\obeylines\newlinechar=`\^^M\pass@rg}
{\obeylines\gdef\pass@rg#1{\writ@line\relax #1^^M\hbox{}^^M}%
\gdef\writ@line#1^^M{\expandafter\toks0\expandafter{\striprel@x #1}%
\edef\next{\the\toks0}\ifx\next\em@rk\let\next=\endgroup\else\ifx\next\empty%
\else\immediate\write\wfile{\the\toks0}\fi\let\next=\writ@line\fi\next\relax}}
\def\striprel@x#1{} \def\em@rk{\hbox{}} {\catcode`\%=12\xdef\pctsign{
\def\semi{;\hfil\break}
\def\addref#1{\immediate\write\rfile{\noexpand\item{}#1}} 
\def\listrefs{\vfill\eject\immediate\closeout\rfile
\baselineskip=12pt\leftline{{\bf References}}\bigskip{\frenchspacing%
\escapechar=` \input refs.tmp\vfill\eject}\nonfrenchspacing}
\def\startrefs#1{\immediate\openout\rfile=refs.tmp\refno=#1}
%
%
\font\titlerm=cmbx12 scaled\magstep5\font\titlerms=cmbx10 scaled\magstep5
\font\titlermss=cmbx7 scaled\magstep3 \font\titlei=cmmi10 scaled\magstep3
\font\titleis=cmmi7 scaled\magstep4 \font\titleiss=cmmi5 scaled\magstep4
\font\titlesy=cmsy10 scaled\magstep4 \font\titlesys=cmsy7 scaled\magstep4
\font\titlesyss=cmsy5 scaled\magstep4 \font\titleit=cmsl10 scaled\magstep4
\skewchar\titlei='177 \skewchar\titleis='177 \skewchar\titleiss='177
\skewchar\titlesy='60 \skewchar\titlesys='60 \skewchar\titlesyss='60
\font\cap=cmr17 scaled\magstep2\font\caps=cmr12 scaled\magstep2
\def\titlefont{\def\rm{\fam0\titlerm}
\textfont0=\titlerm \scriptfont0=\titlerms \scriptscriptfont0=\titlermss
\textfont1=\titlei \scriptfont1=\titleis \scriptscriptfont1=\titleiss
\textfont2=\titlesy \scriptfont2=\titlesys \scriptscriptfont2=\titlesyss
\textfont\itfam=\titleit \def\it{\fam\itfam\titleit} \rm}
\def\abstractfont{\tenpoint}

\def\tenpoint{\def\rm{\fam0\tenrm}
\textfont0=\tenrm \scriptfont0=\sevenrm \scriptscriptfont0=\fiverm
\textfont1=\teni  \scriptfont1=\seveni  \scriptscriptfont1=\fivei
\textfont2=\tensy \scriptfont2=\sevensy \scriptscriptfont2=\fivesy
\textfont\itfam=\tenit \def\it{\fam\itfam\tenit}
\textfont\bffam=\tenbf \def\bf{\fam\bffam\tenbf} \rm}
%
%
\def\newcap#1#2{
\global\advance\secno by1
\xdef\secsym{\the\secno.}\global\meqno=1
  \vskip 18pt
\noindent
{\secfont  #1.\  #2}
\vskip 3mm}
\xdef\secsym{}
%
	\def\subsect#1#2{\advance\subsecnum by 1\subsubsecnum=0
  			           \vskip 14pt
                 \leftline{\subsecfont \the\secno.#1. #2}
           			\vskip -8pt
                 \message{#1}\nobreak \noindent}
%
\font\arm=cmbx12 scaled\magstep2 \font\arms=cmbx10 scaled\magstep2
\font\armss=cmbx7 scaled\magstep2 \font\ai=cmmi10 scaled\magstep2
\font\ais=cmmi7 scaled\magstep3 \font\aiss=cmmi5 scaled\magstep3
\font\asy=cmsy10 scaled\magstep3 \font\asys=cmsy7 scaled\magstep3
\font\asyss=cmsy5 scaled\magstep3 \font\ait=cmsl10 scaled\magstep3
\skewchar\ai='177 \skewchar\titleis='177 \skewchar\aiss='177
\skewchar\asy='60 \skewchar\titlesys='60 \skewchar\asyss='60
\def\afont{\def\rm{\fam0\arm}
\textfont0=\arm \scriptfont0=\arms \scriptscriptfont0=\armss
\textfont1=\ai \scriptfont1=\ais \scriptscriptfont1=\aiss
\textfont2=\asy \scriptfont2=\asys \scriptscriptfont2=\asyss
\textfont\itfam=\ait \def\it{\fam\itfam\ait} \rm}
%
%
\def\title#1{
\vskip 0.4cm
\noindent
\centerline{\bf #1}
\vskip 0.1cm}
%
%
\font\brm=cmbx12 \font\brms=cmbx10
\font\brmss=cmr10  \font\bi=cmmi8
\font\bis=cmmi5  \font\biss=cmmi5
\font\bsy=cmsy7  \font\bsys=cmsy7
\font\bsyss=cmsy7  \font\bit=cmsl10
\skewchar\bi='177 \skewchar\titleis='177 \skewchar\biss='177
\skewchar\bsy='60 \skewchar\titlesys='60 \skewchar\bsyss='60
\def\bfont{\def\rm{\fam0\brm}
\textfont0=\brm \scriptfont0=\brms \scriptscriptfont0=\brmss
\textfont1=\bi \scriptfont1=\bis \scriptscriptfont1=\biss
\textfont2=\bsy \scriptfont2=\bsys \scriptscriptfont2=\bsyss
\textfont\itfam=\bit \def\it{\fam\itfam\bit} \rm}
%
%

%
\def\subtitle#1{
\global\ftno=0
\xdef\secsym{\the\secno.}
\vskip 1cm
\centerline{\bf\afont #1}
\vskip 0.5cm}
\xdef\secsym{}
\font\srm=cmr8 \font\srms=cmr5
\font\srmss=cmr5
\def\sfont{\def\rm{\fam0\srm}
\textfont0=\srm \scriptfont0=\srms \scriptscriptfont0=\srmss
\rm}
\def\feet#1{
\baselineskip=10pt
{\hskip -10mm\sfont #1}
}
%
%
\def\sifont{\def\it{\fam0\srm}
\textfont0=\srm \scriptfont0=\srms \scriptscriptfont0=\srmss
\it}
\def\feit#1{
{\sifont #1}
}
%
%
\def\parsl{\raise.15ex\hbox{/}\kern-.57em\partial}
\def\pasl{\partial \kern-1.45mm {/}}
\def\ii{\'{\i}}
\def\sc{\scriptscriptstyle}
\def\s{* d_a}
\def\lb{\lbrack}
\def\rb{\rbrack}
\def\O{\cal O}
\def\r{\rho}
\def\e{\eta}
\def\s{\sigma}
\def\G{\gamma}
\def\lr{\lbrace}
\def\rr{\rbrace}
\def\ep{\epsilon}
\def\v{\vert}
\def\t{\theta}
\def\lf{\lfloor}
\def\rf{\rfloor}
\def\l{\lambda}
\def\L{\Lambda}
%
%
\font\ptlerm=cmr17 scaled\magstep4\font\ptlerms=cmr17 scaled\magstep4
%
%
$\,$
\overfullrule=0pt
%
%
\def\be{\beta}
\def\ch{\chi}
\def\ga{\gamma}
\def\de{\delta}
\def\ep{\varepsilon}
\def\ze{\zeta}
\def\io{\iota}
\def\ka{\kappa}
\def\la{\lambda}
\def\na{\nabla}
\def\ro{\varrho}
\def\si{\sigma}
\def\om{\omega}
\def\ph{\varphi}
\def\ta{\tau}
\def\th{\theta}
\def\te{\vartheta}
\def\up{\upsilon}
\def\Ga{\Gamma}
\def\De{\Delta}
\def\La{\Lambda}
\def\Si{\Sigma}
\def\Om{\Omega}
\def\Te{\Theta}
\def\Th{\Theta}
\def\Up{\Upsilon}
\def\CO{{\cal{O}}}
\def\sign{{\rm sign}}
\def\inbar{\,\vrule height1.5ex width.4pt depth0pt}
\def\IC{\relax{\raise .7mm\hbox{${\scriptstyle\vert}$}}\hskip -1.5mm {\rm C}}
\def\IN{\relax{\raise .4mm\hbox{${\scriptscriptstyle\vert}$}}\hskip 
-1.4mm {\rm N}}

\def\R{\scriptscriptstyle R}
\def\scc{\scriptscriptstyle}
\def\ba{${\scriptstyle /}$}
\def\ve{${\scriptstyle \vert}$}
\def\bas{${\scriptscriptstyle /}$}
\def\Cs{C \kern-2.2mm \raise.4ex\hbox {\ba} \kern.2mm}
\def\Css{C \kern-1.9mm \raise.3ex\hbox {\bas} \kern.8mm}
\def\Zs{Z \kern-2.2mm \raise.4ex\hbox {\ba} \kern.2mm}
\def\Rs{{\rm R} \kern-2.9mm \raise.35ex\hbox {\ve} \kern.4mm\hskip 1.3mm}
\def\Ns{{\rm N} \kern-2.9mm \raise.35ex\hbox {\ve} \kern.4mm\hskip 1.3mm}
\def\Qs{{\rm Q} \kern-2.4mm \raise.35ex\hbox {\ve} \kern.4mm\hskip 1.3mm}
\def\boxe{\sqcup \kern-2.3mm \sqcap}
%


\hsize = 14.6truecm
\vsize = 21.0truecm

\nopagenumbers

\headline{
            \rm \hfill\the\pageno}

\input epsf


\hoffset=6truemm
\baselineskip=13pt
\parskip 2pt plus 1pt
\parindent=15pt

\font\notefont=cmr9

\vbadness=10000
\widowpenalty=10000
\clubpenalty=10000

\font\afont=cmbx12
\font\bfont=cmr9
\font\cfont=cmsl9
\font\hfont=cmcsc10
\font\pfont=cmss10
\font\tfont=cmbx12
\font\reffont=cmcsc8
\font\secfont=cmbx10
\def\subsecfont{\bf}
\def\subsubsecfont{\sl}
\font \abscwifont=cmss10

\font\footfont=cmr8
\def\ref#1{{\reffont#1}}

\def\sq{\quad}
\def\sqq{\qquad}
\def\cl{\centerline}

\newcount\secnum
\newcount\subsecnum
\newcount\subsubsecnum
\newcount\remarknum
\newcount\lemmanum
\newcount\thmnum
\newcount\vbnum

%
%
%

\def\subsubsect#1{\advance\subsubsecnum by 1
							\vskip 12pt
 
\leftline{\subsubsecfont\the\secnum.\the\subsecnum.\the\subsubsecnum. 
\ #1}
           			\vskip -6pt
                \message{#1}\nobreak \noindent}

\def\remark{\advance \remarknum by 1
          \vskip\baselineskip
          \noindent {\bf Remark \the\secnum.\the\remarknum.}\quad\rm
          }
\def\lemma{\advance \lemmanum by 1
          \vskip\baselineskip
          \noindent {\bf Lemma \the\secnum.\the\lemmanum.}\quad\it
          }
\def\theorem{\advance \thmnum by 1
          \vskip\baselineskip
          \noindent {\bf Theorem \the\secnum.\the\thmnum.}\quad\it
          }
\def\proof{
          \vskip\baselineskip
          \noindent {\bfProof }\quad\rm
          }
\def\vb{\advance \vbnum by 1
          \vskip\baselineskip
          \noindent {\bf Example 
\the\vbnum.}\quad
          }

\def\tfrac#1#2{{{\lower.4ex
\hbox{$\scriptstyle#1$}}\over
{\raise.4ex
\hbox{$\scriptstyle#2$}}}}

\def\remp
#1. #2\par{\medbreak\noindent{\bf#1.\enspace}{\rm#2}\par\medbreak}

\def\thep
#1. #2\par{\medbreak\noindent{\bf#1.\enspace}{\sl#2}\par\medbreak}

\def\prop
#1. #2\par{\medbreak\noindent{\bf#1.\enspace}{\rm#2}\par
  \rightline{\vrule height4pt width5.5pt depth2pt}\medbreak}
\def\proof{\bf \medbreak \noindent Proof. \rm}
\def\eoproof{{\unskip\nobreak\hfil\penalty50
   \hskip2em\hbox{}\nobreak\hfil\vrule height4pt width5.5pt depth2pt
   \parfillskip=0pt\finalhyphendemerits=0\medbreak}}

\def\w#1{{\sqrt#1}\,}
\def\kd{\partial}
\def\sq{\quad}
\def\sqs{\qquad}
\def\n{\eqno}
\def\cl{\centerline}
\def\el{\eqalign}

\def\iy{\infty}

\def\bo{{\cal O}}

\def\pd#1#2{{{\kd#1}\over{\kd#2}}}
\def\pdt#1#2{{{\kd^2#1}\over{\kd #2^2}}}

\def\frac#1#2{{{#1}\over{#2}}}

\def\intp{\int_0^\iy}
\def\intr{\int_{-\iy}^\iy}

\def\RR{{{\rm I}\!{\rm R}}}
\def\NN{{{\rm I}\!{\rm N}}}
\def\RRP{{{\rm I}\!{\rm R^+}}}
\def\RRN{{{\rm I}\!{\rm R_0}}}
\def\NNP{{{\rm I}\!{\rm N^+}}}
\def\ZZ{{\hbox{Z}\!\!\hbox{Z}}}
\def\KK{{{\rm I}\!{\rm K}}}
\def\one{{{\rm 1}\hskip-0.55ex{\rm I}}}
\def\CC{\hbox{\rlap{$\,\,
   $\hbox{\vrule height6.2pt width.35pt depth-0.1pt}}$\rm C$}}
\def\QQ{\hbox{\rlap{$\,\,
   $\hbox{\vrule height6pt width.35pt depth0.1pt}}$\rm Q$}}
\def\P{\cal P}

%
\newcount\equationnumber	\equationnumber=0
\def\eqnum{\relax
	\global\advance\equationnumber by 1
	\equationnumberformat{\the\equationnumber}%
	}%
\def\eqname#1{\relax
	\count255=\equationnumber
	\assignnumber{EN#1}\equationnumber
	\global\equationnumber=\count255
	\global\advance\equationnumber by 1
	\ifnum\csname EN#1\endcsname=\equationnumber
	\else
		\message{The equation number for ``#1'' is incorrect!}%
	\fi
	\equationnumberformat{\csname EN#1\endcsname}%
	}%
\def\equationnumberformat#1{\eqno(\the\secnum.\equationnumbertype{#1})}%
\def\equationnumbertype#1{\number#1\relax}%
\def\referenceequation#1{\relax
	\assignnumber{EN#1}\equationnumber
	\equationnumbertype{\csname EN#1\endcsname}%
	}%
\def\forwardreferenceequation#1#2{\relax
	\global\advance\equationnumber by #2
	\assignnumber{EN#1}\equationnumber
	\global\advance\equationnumber by -1
	\global\advance\equationnumber by -#2
	\referenceequation{#1}%
	}%
%
\def\assignnumber#1#2{\relax
	\ifnum0<0\csname#1\endcsname
	\else
		\global\advance#2 by 1
		\expandafter\expandafter\expandafter
			\xdef\csname#1\endcsname{\the#2}%
	\fi
	}%

\def\fre{\forwardreferenceequation}
\def\en{\eqname}
\def\req#1{(\the\secnum.{\referenceequation{#1}})}

\def\arcsinh{{\rm arcsinh}}
\def\arccosh{{\rm arccosh}}
\def\arctanh{{\rm arctanh}}
\def\arccoth{{\rm arccoth}}

\def\erf{{\rm erf}}
\def\erfc{{\rm erfc}}

\def\phase{{\rm phase}}
\def\sign{{\rm sign}}

\def\wt{{\sqrt{2}}}

\def\phih{{\widehat \phi}}
\def\psih{{\widehat \psi}}

\def\Ai{{{\rm Ai}}}
\def\Bi{{{\rm Bi}}}

\def\phizeta{\left(\frac{4\z}{1-z^2}\right)^{1/4}}

\catcode`@=11 

\font\ninerm=cmr10 at 9pt
\font\eightrm=cmr7 at 8pt
\font\sixrm=cmr7 at 6pt

\font\ninei=cmmi10 at 9pt
\font\eighti=cmmi10 at 8pt
\font\sixi=cmmi7 at 6pt
\skewchar\ninei='177
\skewchar\eighti='177
\skewchar\sixi='177

\font\ninesy=cmsy10 at 9pt
\font\eightsy=cmsy10 at 8pt
\font\sixsy=cmsy7 at 6pt
\skewchar\ninesy='60
\skewchar\eightsy='60
\skewchar\sixsy='60

\font\eightss=cmssq8

\font\eightssi=cmssqi8

\font\ninebf=cmbx10 at 9pt
\font\eightbf=cmbx10 at 8pt
\font\sixbf=cmbx7 at 6pt

\font\ninett=cmtt10 at 9pt
\font\eighttt=cmtt10 at 8pt

\hyphenchar\tentt=-1 
\hyphenchar\ninett=-1
\hyphenchar\eighttt=-1

\font\ninesl=cmsl10 at 9pt
\font\eightsl=cmsl10 at 8pt

\font\nineit=cmti10 at 9pt
\font\eightit=cmti10 at 8pt

\font\tenu=cmu10 
\newskip\ttglue

\def\ninepoint{\def\rm{\fam0\ninerm}%
  \textfont0=\ninerm \scriptfont0=\sixrm \scriptscriptfont0=\fiverm
  \textfont1=\ninei \scriptfont1=\sixi \scriptscriptfont1=\fivei
  \textfont2=\ninesy \scriptfont2=\sixsy \scriptscriptfont2=\fivesy
  \textfont3=\tenex \scriptfont3=\tenex \scriptscriptfont3=\tenex
  \def\it{\fam\itfam\nineit}%
  \textfont\itfam=\nineit
  \def\sl{\fam\slfam\ninesl}%
  \textfont\slfam=\ninesl
  \def\bf{\fam\bffam\ninebf}%
  \textfont\bffam=\ninebf \scriptfont\bffam=\sixbf
   \scriptscriptfont\bffam=\fivebf
  \def\tt{\fam\ttfam\ninett}%
  \textfont\ttfam=\ninett
  \tt \ttglue=.5em plus.25em minus.15em
  \normalbaselineskip=11pt
  \def\MF{{\manual hijk}\-{\manual lmnj}}%
  \let\sc=\sevenrm
  \let\big=\ninebig
  \setbox\strutbox=\hbox{\vrule height8pt depth3pt width\z@}%
  \normalbaselines\rm}

\def\eightpoint{\def\rm{\fam0\eightrm}%
   \textfont0=\eightrm \scriptfont0=\sixrm \scriptscriptfont0=\fiverm
   \textfont1=\eighti \scriptfont1=\sixi \scriptscriptfont1=\fivei
   \textfont2=\eightsy \scriptfont2=\sixsy \scriptscriptfont2=\fivesy
   \textfont3=\tenex \scriptfont3=\tenex \scriptscriptfont3=\tenex
   \def\it{\fam\itfam\eightit}%
   \textfont\itfam=\eightit
   \def\sl{\fam\slfam\eightsl}%
   \textfont\slfam=\eightsl
   \def\bf{\fam\bffam\eightbf}%
   \textfont\bffam=\eightbf \scriptfont\bffam=\sixbf
    \scriptscriptfont\bffam=\fivebf
   \def\tt{\fam\ttfam\eighttt}%
   \textfont\ttfam=\eighttt
   \tt \ttglue=.5em plus.25em minus.15em
   \normalbaselineskip=9pt
   \def\MF{{\manual opqr}\-{\manual stuq}}%
   \let\sc=\sixrm
   \let\big=\eightbig
   \setbox\strutbox=\hbox{\vrule height7pt depth2pt width\z@}%
   \normalbaselines\rm}

\def\b{\beta}\def\c{\gamma}\def\d{\delta} \def\vth{\vartheta}
\def\eps{\varepsilon}\def\f{\phi}\def\k{\kappa}\def\l{\lambda}\def\m{\mu}
\def\p{\pi}\def\r{\rho}\def\s{\sigma}\def\t{\tau}\def\th{\theta}
\def\x{\xi}\def\y{\eta}\def\z{\zeta}\def\om{\omega}
\def\La{\Lambda}\def\oom{\Omega}\def\G{\Gamma}\def\D{\Delta}

\def\sn{\sum_{n=0}^\iy}
\def\sk{\sum_{k=0}^\iy}
\def \som#1#2{\sum_{{#1}}^{{#2}}}

\def\nfe{\eject}
\def\vfe{\vfil\eject}
\def\oom{\Omega}
\def\wh{\widehat}

\def\gp{$C_n^{\c}(x)$}
\def\lp{$L_n^{\a}(x)$}
\def\jp{$P_n^{(\a,\b)}(x)$}

\def\gpf{C_n^{\c}(x)}
\def\lpf{L_n^{\a}(x)}
\def\jpf{P_n^{(\a,\b)}(x)}

\def\bri{\frac1{2\pi i}}
\def\C{{\cal C}}
\def\H{{\cal H}}
\def\L{{\cal L}}

\def\({\left(}
\def\){\right)}

\def\[{\left[}
\def\]{\right]}

\def\aa{^{(\a)}(x)}

\noindent
%

%% file: 102153ref.tex
\nref\bowong{{\reffont Rui Bo and R. Wong},
Uniform asymptotic expansion of Charlier polynomials.
{\it Methods Appl. Anal.} {\bf 1} (1994), 294-313.}

\nref\dunster{{\reffont T. M. Dunster},
Uniform asymptotic expansion for Charlier polynomials.
{\it J. Approx. Theor.} {\bf 112} (1994), 93-133.}

\nref\bateman{{\reffont A. Erd\'elyi},
{\it Higher transcendental functions},
Vols. I, II, III, A. Erd\'elyi, ed.
{ McGraw-Hill}
(Reprinted by Krieger Publishing Co., Malabar, FL,  1981.) (1953).}

\nref\frong{{\reffont C.L. Frenzen and R. Wong},
Uniform asymptotic expansions of Laguerre polynomials,
{\it SIAM J. Math. Anal.}, {\bf 19} (1988),  1232--1248.}

\nref\goh{{\reffont W.M.Y. Goh},
Plancherel-Rotach asymptotics for the Charlier polynomials.
{\it Constructive Approximation} {\bf 14}(1998), 151-168.}

\nref\hsu{{\reffont L.C. Hsu},
Certain asymptotic expansions for Laguerre polynomials
and Charlier polynomials.
{\it Approx. Theory Appl. (N.S.)} {\bf 11} (1995), 94-104. }

\nref\koek{{\reffont R. Koekoek and R. F. Swarttouw},
{\it Askey scheme of hypergeometric orthogonal polynomials and its $q-$analogue},
{\tt http://aw.twi.tudelft.nl/\~{}koekoek/askey.html}, (1999).}

\nref\temlopherm{{\reffont J. L. Lopez and N.M. Temme},
Approximations of orthogonal polynomials
in terms of Hermite polynomials.
{\it Methods Appl. Anal.} {\bf 6} (1999), 131-146.}

\nref\temmelopez{{\reffont J. L. Lopez and N.M. Temme},
Two-point Taylor expansions of analytic functions.
{\it Stud. Appl. Math}, {\bf 109}  (2002), 297-311.}

\nref\multayl{{\reffont Jos\'e L. Lopez and Nico M. Temme},
Multi-point Taylor expansions of analytic functions.
Submitted for publication.}

\nref\luke{{\reffont Y. L. Luke}, 
{\it The Special Functions and their
Applications}, {Academic Press, New York} (1969).}

\nref\temlag{{\reffont N.M. Temme},
Asymptotic estimates for Laguerre polynomials,
{\it ZAMP} {\bf 41} (1990), 114--126.}

\nref\temsf{{\reffont N.M. Temme},
{\it Special Functions. An Introduction to the
Classical Functions of Mathematical Physics},
{John Wiley and Sons, New York} (1996).}

\nref\tric{{\reffont F.G. Tricomi},
{\it Funzioni ipergeometriche confluenti},
{Edizione Cremonese, Roma} (1954)}.

\nref\wong{{\reffont R. Wong}, 
{\it Asymptotic Approximations of
Integrals}, {Academic Press, New York}  (1989).
Reprinted by SIAM in 2001.}